\input{ALT-TITELEI/enskog.top}
\documentclass[a4paper, 
               10pt, 
               oneside 
               ]{article} 
\usepackage{
            amsmath,     
            amsfonts,    
            amsthm,      
            textcomp,    
            graphics,    
            epsfig,      
            pstricks,      
            makeidx,     
            ifthen,
            verbatim,
            latexsym,
            exscale      
           }
\usepackage[dvips, 
            colorlinks,
            pagebackref,
            hyperindex,
]{hyperref}

\catcode`\@=11
  \def\ps@headings{%
    \let\@oddfoot\footline   
    \def\@oddhead{\hfil\thepage}
    \let\@mkboth\markboth
  }
\catcode`\@=12
\catcode`\@=11
\DE{
}
\EN{
}
\catcode`\@=12

\newif\ifamstheorem 
\amstheoremtrue
\newif\ifsptheorem 
\newif\ifnochaps
\nochapstrue
\usepackage{ALT-TITELEI/LATEX/refsect} %
\usepackage{ALT-TITELEI/LATEX/layout} 
\usepackage{ALT-TITELEI/LATEX/math} 
\usepackage{ALT-TITELEI/LATEX/mathdisplay} %
\usepackage{pifont} 
\def\actcolor{black} 
\ifnocolor\renewcommand{\color}[1]{}\fi
\usepackage{ALT-TITELEI/LATEX/msymbol} 
\makeindex

\usepackage{ALT-TITELEI/script}
\usepackage{ALT-TITELEI/enskog}

\begin{document}

{
\pagestyle{empty}

\begin{center}

\vspace*{13mm}\noindent
{\Large\bf Relativistic equations for the}
\\[1.2ex]
{\Large\bf generalized Chapman-Enskog hierarchy.} 
\\
\vspace*{10mm}\noindent
{\large\bf Hans Wilhelm Alt}
\\
{\bf Technische Universit\"at M\"unchen}
\\
\vspace*{10ex}\noindent
The paper has been accepted by AMSA
(Advances in Mathematical Sciences and Applications) 
and will probably be 
published in Vol. 25 (2016).
\\
\end{center}

}

\vspace*{15ex}\noindent


\setcounter{tocdepth}{2} 
\catcode`\@=11
\@starttoc{toc}
\catcode`\@=12
\newpage
{}
\sect{intro}{\DE{}\EN{Introduction}}

In the book \CIT{ChapmanCowling}{}  
Chapman \& Cowling have presented a theory based
on Boltzmann's equation of monatomic gases.
Some of the results one finds in section \REF{chapman..}.
The Chapman-Enskog system is a hierarchy of conservation laws
for the mass, the momentum, the second moments, the third moments
and so on. The hierarchy of functions is denoted by $F_{i_1\cdots i_{N+1}}$
and the conservation laws are as one can see
in \CIT{Mueller1998}{Chap. 2 (3.15)} or in this paper
in equation \EQU{higher.*.sys-coriolis}
\begin{Equ}{sys-coriolis}
  \sum_{j\geq0}\de_{y_j}F_{\alpha j}
  -\sum_{\beta\in\{0\til3\}^{N+1}}\coriolis^\beta_\alpha F_\beta=\tx\force_\alpha
  \quad\tm{for} \alpha\in\{0\til3\}^N
\end{Equ}
where $\beta=i_1\cdots i_{N+1}$.
It is important to mention, that in this paper, except section \REF{chapman..},
the higher moments $F_{i_1\cdots i_{N+1}}$
are defined as independent variables, so that, for example, the pressure $p$ and the internal energy $\eps$ are independent
from each other and not as in \EQU{chapman.*.special}.
Only a constitutive relation between them based on Gibb's relation
is assumed, therefore we have the general form
\[\begin{arr}c
  F_{k0}=\rho\tx v_k \,,\quad F_k=\rho\tx v_k+\JJ_k
\,,\\
  F_{kl0}=\rho\tx v_k\tx v_l+\eps_{kl} \,,\quad
  F_{kl}=\rho\tx v_k\tx v_l+\Pi_{kl}
\,,\end{arr}\]
which includes Navier-Stokes equations.
Also the higher equations with terms for $N\geq2$ are in the same way a
generalization of the Navier-Stokes equation.
\DE{}
\EN{The possible symmetry of $F_{i_1\cdots i_{N+1}}$
only refers to the first $N$ indices.}

\medskip
The higher moments $F_{i_1\cdots i_M}$, here $M=N+1$, satisfy an identity for different observers
which is crucial for the entire theory. This identity says that $F_{i_1\cdots i_M}$ are contravariant $M$-tensors (see definition \EQU{contra})
\begin{Equ}{rule-F}
  F_{i_1\cdots i_M}\circ Y
  =\sum_{\bar i_1\til\bar i_M\geq0}
  Y_{i_1\p\bar i_1}\cdots Y_{i_M\p\bar i_M}
  F^*_{\bar i_1\cdots\bar i_M}
\end{Equ}
as in \EQU{chapman.trans.o} and \EQU{higher.*.rule-T}.
This means that one has an observer transformation
between two observers $y=Y(y^*)$, where,
to have a common description,
$y\in\RR^4$ are the time and space coordinates
of one observer and $y^*\in\RR^4$ the coordinates of a second observer,
and finally $Y$ denotes the observer transformation.
Here $y=(t,x)\in\RR^4$ in the classical sense.
Now, the moments $F$ are the quantities of one observer 
and the same moments $F^*$ the quantities for the other observer. 
Similar is the notion for other quantities like the forces $\tx\force_\alpha$
which are denoted by $\tx\force_\alpha^*$ for the other observer.
It is important that
the differential equations of the system are the same
for all observers, that is, the system \EQU{sys-coriolis} for
$F_\beta$, $\coriolis^\beta_\alpha$, $\tx\force_\alpha$ is the same for
$F^*_\beta$, $\coriolis^{*\beta}_\alpha$, $\tx\force^*_\alpha$,
see section \REF{scalar..} and section \REF{higher..}.

\medskip
The difference between the classical formulation and the relativistic version
lies in the Group of transformations, which in the relativistic case
is based on Lorentz matrices (see section \REF{invariance..}).
It is also based on a matrix $\GG$ which occurs in the conservation laws
and is transformed by the rule $\GG\circ Y=\D{Y}\GG^*\up{(\D{Y})}T$,
that is, $\GG$ is a contravariant tensor.
There are special cases \[
  \GG=\GG_\cc=\begin{mat}-\frac1{\cc^2}\tab 0\\0\tab\Id\end{mat}
\quad\tm{and}\quad
  \GG=\GG_\infty=\begin{mat}0\tab 0\\0\tab\Id\end{mat}
,\]
where $\Id$ is the Identity in space
and where $\GG_\infty$ is the limit of $\GG_\cc$ as $\cc\to\infty$,
hence $\GG_\infty$ is the matrix in the classical limit.
Therefore the matrix
$\GG_\cc$ in the relativistic case is invertible
whereas $\GG_\infty$ is not. This limit is justified although the speed of light
in vacuum
\[\cc=2.99792458\cdot10^8\frac{m}{s}\] is a finite number.
The inverse matrix $\up{\GG}{-1}$, if it exists, has the transformation rule 
$\up{\GG^*}{-1}=\up{(\D{Y})}T\up{\GG}{-1}\circ Y\D{Y}$,
that is, $\up{\GG}{-1}$ is a covariant tensor.
And $\up{\GG}{-1}$ is the commonly used matrix in relativistic physics.

\medskip
It is an essential step that
in section \REF{time..} we introduce a covariant vector $\ewelt$
with the identity
\[
  \ewelt\dd\GG\ewelt=-\frac1{\cc^2}
\,.\]
This vector is the ``time direction'' and plays an important role.
With $\ewelt=e'_0$ it is part of a dual basis
$\{e'_0,e'_1,e'_2,e'_3\}$ of $\{e_0,e_1,e_2,e_3\}$, and
the matrix $\GG$ has the non-unique representation \REF{time.prop.G}
\begin{Equ}{prop-G}
  \GG=-\frac1{\cc^2} e_0\up{e_0}T+\sum_{i=1}^ne_i\up{e_i}T
\,.\end{Equ}
Therefore $\ewelt$ is strongly connected to the matrix $\GG$
and because $\ewelt^*=\up{\D{Y}}T\ewelt\circ Y$ it can occur in
basis physical laws. One fundamental example is the definition
of a 4-velocity $\tx v$ with $\ewelt\dd\tx v=1$ in \EQU{scalar.def-velo.o}.
This definition differs by a scalar mul\-ti\-plication from the
known definition (see \REF{scalar.u.}).
Our definition of a 4-velocity plays an important role
in the constitutive law \REF{fluid.Dv.} for fluids.

\medskip
There is an essential application of the time vector $\ewelt$,
it is the reduction principle.
It is obseved in \REF{chapman.reduct.} that
in the classical limit the system of $N^{\rm th}$-order moments contain
the system for $(N-1)^{\rm th}$-order moments.
This classical reduction principle is generalized
in this paper to the relativistic case.
In section \REF{impuls..} we treat the case $N=1$, that is,
the 4-momentum equation.
Here the simplest case of a reduction occurs,
choosing $\zeta=\eta\ewelt$ as test function we are able to show
that the mass equation is part of the 4-momentum equation,
see \REF{impuls.mass.}.
This prove in the general version is then presented to the equation of $N^{\rm th}$-moments
in \REF{higher.reduction.}. We use as test function
$\zeta_{\alpha_1\cdots\alpha_N}:=\ewelt_{\alpha_1}\eta_{\alpha_2\cdots\alpha_N}$
and obtain this way a system of $(N-1)^{\rm th}$-moments. Therefore the relativistic version of the $N^{\rm th}$-order moments
is a generalization of the classical case.

\medskip
One main tool in this paper is the principle of relativity.
Its application to basic differential equations is presented
in section \REF{continuum..}. The quantities of these
conservation laws have to satisfy certain transformation rules
\EQU{continuum.*.o}.
This is applied to scalar conservation laws in \EQU{scalar.*.o},
to the 4-momentum equation in \EQU{impuls.*.rule}
and to $N^{\rm th}$-moment equations in 
\EQU{higher.*.rule-T} and \EQU{higher.*.rule-g}.
We mention that for Maxwell's equation in vacuum
in \CIT{Einstein1905}{II Elektrodynamischer Teil}
this method
has been applied to Lorentz transformations.

\medskip
In section \REF{point..} we require that the mass-momentum equation
of section \REF{impuls..} should also describe the law of particles.
Thus one has to consider the differential equation
$\tx\div\,T = \tx r$ for distributions, here a
one-dimensional curve $\Gamma$ in spacetime $\RR^4$.
So $\Gamma:=\set{\tx\xi(s)}{s\in\RR}$ is the evolving point,
by which we mean that $s$ is chosen so that
$\ewelt(\tx\xi(s))\dd\de_s\tx\xi(s)>0$.
It is shown in a rigorous way
that the differential equation is equivalent to known ordinary
differential equations.

\medskip
As a concluding contribution
let us write down the definitions of a contravariant $m$-tensor
$T=\seq{T_{k_1\cdots k_m}}{k_1\til k_m}$ 
\begin{Equ}{contra}
   T_{k_1\cdots k_m}\circ Y=\sum_{\bar k_1\til\bar k_m\geq0}
  Y_{k_1\p\bar k_1}\cdots Y_{k_m\p\bar k_m}T^*_{\bar k_1\cdots\bar k_m}
\,,\end{Equ}
and the definition of a covariant $m$-tensor
$T=\seq{T_{k_1\cdots k_m}}{k_1\til k_m}$ 
\begin{Equ}{co}
  T^*_{\bar k_1\cdots\bar k_m}=\sum_{k_1\til k_m\geq0}
  Y_{k_1\p\bar k_1}\cdots Y_{k_m\p\bar k_m} T_{k_1\cdots k_m}\circ Y
\,.\end{Equ}
In this connection a 4-matrix is a $2$-tensor
and a 4-vector a $1$-tensor. Besides this we 
call a scalar quantity $u$ ``objective'' if 
$u\circ Y=u^*$ is true.
We denote with an underscore terms in spacetime
which are usually meant in space only, so
$\tx\div q=\sum_{i\geq0}\de_{y_i}q_i$ which in coordinates $y=(t,x)$
means $\tx\div q=\de_tq_0+\sum_{i\geq1}\de_{x_i}q_i$.
Also $\tx v(y)\in\RR^4$ is the spacetime version of the velocity. 

{}
{}
\sect{chapman}{Chapman-Enskog \DE{}\EN{method}}
\newcommand{\spforce}{\mathbf{g}}
\newcommand{\fc}{\mathbf{g}}
\DE{}
\EN{The classical Boltzmann equation is a differential equation
for the probability $(t,x,c)\mapsto f(t,x,c)$.
For a single species this equation reads as}
\begin{Equ}{boltzmann}\begin{arr}c
  \de_tf+\sum_{i=1}^3c_i\de_{x_i}f+\sum_{i=1}^3\spforce_i\de_{c_i}f
  =\rate
\end{arr}\end{Equ} 
with the additional equation $\sum_{i=1}^n\de_{c_i}\spforce_i=0$
for the external acceleration $\spforce$.
The quantity $f$ is the density of atoms at $(t,x)$ with velocity $c$,
and the acceleration $\spforce$ is a function of $(t,x,c)$.
Moreover, $\rate$ is the collision product,
which is also a function of $(t,x,c)$.
The Boltzmann equation is
explained in many papers including the collision product,
see for example \CIT{Mueller1985}{5.2.1}
and the literature cited there.

\medskip
\DE{}\EN{According to the probability $f$ the \DEF{higher moments}
are defined for $k_1\til k_M\in\{0\til 3\}$ by}
\begin{Equ}{moments-def}
  F_{k_1\cdots k_M}(t,x)
  :=\int_{\RR^3}m\tx c_{k_1}\cdots \tx c_{k_M}f(t,x,c)\d{c}
\,,\end{Equ}
\DE{}
\EN{where $m$ is the particle mass and $\tx c$ the extended velocity}
\[
  \tx c:=\begin{mat}1\\c\end{mat}
\,,\]
\DE{}\EN{that is}
$\tx c=(c_0,c)=(c_0,c_1\til c_3)$, $c_0:=1$.
\DE{}
\EN{We remark that $F_{k_1\cdots k_M}$ can be the same function for different
indices, for example $F_{k0}=F_k$ and $F_{0kl}=F_{kl}$.}
\DE{}
\EN{If $f$ is a solution of Boltzmann's equation
and decays fast enough for $|c|\to\infty$
then the higher moments satisfy the following system
of partial differential equations in $(t,x)$:}
\DE{}
\EN{For $i_1\til i_N\in\{0\til 3\}$ there holds}
\begin{Equ}{moments-pde}
  \de_tF_{i_1\cdots i_N}+\sum_{i=1}^3\de_{x_i}F_{i_1\cdots i_Ni}
  =R_{i_1\cdots i_N}
\,,\end{Equ}
\begin{Equ}{moments-rs}\begin{arr}l
  R_{i_1\cdots i_N}(t,x)
  :=\int_{\RR^3}m\tx c_{i_1}\cdots \tx c_{i_N}\rate(t,x,c)\d{c}
\\\hfill
  +\sum_{i}\int_{\RR^3}m\spforce_i(t,x,c)
  \de_{c_i}(\tx c_{i_1}\cdots \tx c_{i_N}) f(t,x,c)\d{c}
\end{arr}\end{Equ}
This has been proved in many books, see for example \CIT{Mueller1985}{5.2.2}
and \CIT{Mueller1998}{Chap.2 3.1 3.2 (3.15)}.
The equations for $N\leq2$ are the classical equations,
where $\rho$ is assumed to be positive
and for $N=2$ one uses only the trace:
\begin{Equ}{special}\begin{arr}l
  \rho:=\int_{\RR^3}mf\d{c}=F_0
\,,\quad
  \rho v:=\int_{\RR^3}mfc\d{c}=\seq{F_i}{i=1\til n}
\,,\\
  \Pi:=\int_{\RR^3}mf(c-v)\up{(c-v)}T\d{c}
\,,\quad
  \force:=\int_{\RR^3}mf\spforce\d{c}
\,,\\
  e=\eps+\frac\rho2|v|^2=\frac12\sum_{i=1}^3F_{ii}
  =\int_{\RR^3}\frac{m}2|c|^2f\d{c}
\,,\\
  \eps:=\int_{\RR^3}\frac{m}2|c-v|^2f\d{c}=\frac12\sum_{k=1}^3\Pi_{kk}
\,,\\
  q:=\int_{\RR^3}mf|c-v|^2(c-v)\d{c}
\,,\quad
  g=\int_{\RR^3}m\spforce\dd (c-v) f\d{c}
\,,\end{arr}\end{Equ}
and these quantities satisfy the following classical equations
\begin{equ}{o}\begin{arr}l
  \de_t\rho+\div (\rho v) = 0 
\,,\\
  \de_t(\rho v)+\div (\rho v\up{v}T+\Pi) = \force 
\,,\\
  \de_te+\div (e v + \up{\Pi}{T}v + q) = v\dd\force + g 
\,.\end{arr}\end{equ}
\DE{}
\EN{Here the special properties of the collision term are used,
that is, it does not give any contribution to the mass, momentum, and energy.}
It is clear that the mass equation, which means $N=0$,
is part of the mass-momentum equation, which are the moments with $N=1$.
Similarly one can consider the moments for $N\leq3$ with a trace for $N=3$.
One obtains Grad's 13-moment theory, see \CIT{Grad}{}
and \CIT{Mueller1998}{Chap.2 (3.16) 3.4}.
Then the equations with $N=1$ are contained in the larger system with $N=2$.
\begin{stmt}{reduct}{Observation}
The system of moments \EQU{*.moments-pde} for $i_1\til i_N\in\{0\til 3\}$ 
contains the system for $(N-1)$-order moments.
\end{stmt}
We will not discuss this in detail
but rather focus on the following fact
which makes this observation obvious.
\DE{}
\EN{By denoting the coordinates $y=(t,x)$ the
differential equation \EQU{moments-pde} reads}
\begin{Equ}{moments-dgl}
  \sum_{k=0}^3\de_{y_k}F_{i_1\cdots i_Nk} = R_{i_1\cdots i_N}
\end{Equ}
\DE{}\EN{for} $i_1\til i_N\in\{0\til 3\}$.
\DE{}
\EN{These are all differential equations for moments less or equal $N$.
For these moments the following transformation rule hold,
where the indices run from $0$ and where $M=N+1$.}

\begin{stmt}{trans}{\DE{}\EN{Transformation rule}}
\DE{}\EN{For} $k_1\til k_M\in\{0\til 3\}$
\begin{Equ}{o}
  F_{k_1\cdots k_M}\circ Y
  =\sum_{\bar k_1\til\bar k_M=0}^3
  Y_{k_1\p\bar k_1}\cdots Y_{k_M\p\bar k_M}
  F^*_{\bar k_1\cdots\bar k_M}
\,.\end{Equ}
\end{stmt}
\begin{prf}{}
\DE{}
\EN{For $f$ we have the transformation rule}
$f(t,x,c)=f^*(t^*,x^*,c^*)$
\DE{}\EN{if}
\[
  \begin{mat}t\\x\\c\end{mat}
  =\begin{mat}T(t^*)\\X(t^*,x^*)\\\dot X(t^*,x^*) + Q(t^*)c^*\end{mat}
\,,\]
\DE{}
\EN{where $T(t^*)=t^*+\tdiff$ and $X(t^*,x^*)=Q(t^*)x^*+\xdiff(t^*)$.
Hence with}
\[
  \begin{mat}t\\x\end{mat}
  =Y\left(\begin{mat}t^*\\x^*\end{mat}\right)
  =\begin{mat}T(t^*)\\X(t^*,x^*)\end{mat}
\,,\quad
  \D{Y}=\begin{mat}1\tab0\\\dot X\tab Q\end{mat}
\,,\]
\DE{}
\EN{we obtain the following rule for $F_{k_1\til k_N}$}
\[\begin{arr}l
  F_{k_1\cdots k_N}(t,x)
  =\int_{\RR^3}m\tx c_{k_1}\cdots \tx c_{k_N}f(t,x,c)\d{c}
\\
  =\int_{\RR^3}m\tx c_{k_1}\cdots \tx c_{k_N}f^*(t^*,x^*,c^*)\d{c^*}
\\
  =\int_{\RR^3}m\Big(\sum_{\bar k_1=0}^3Y_{k_1\p\bar k_1}\tx c^*_{\bar k_1}\Big)
  \cdots\Big(\sum_{\bar k_N=0}^3Y_{k_N\p\bar k_N}\tx c^*_{\bar k_N}\Big)
  f^*\d{c^*}
\\
  =\sum_{\bar k_1\til\bar k_N=0}^3
  Y_{k_1\p\bar k_1}\cdots Y_{k_N\p\bar k_N}
  \int_{\RR^3}m\tx c^*_{\bar k_1}\cdots\tx c^*_{\bar k_N}f^*\d{c^*}
\\
  =\sum_{\bar k_1\til\bar k_N=0}^3
  (Y_{k_1\p\bar k_1}\cdots Y_{k_N\p\bar k_N})(t^*,x^*)
  F^*_{\bar k_1\cdots\bar k_N}(t^*,x^*)
\,,\end{arr}\]
\DE{}\EN{where}
\[
  \tx c=\begin{mat}1\\c\end{mat}
  =\begin{mat}1\\\dot X + Qc^*\end{mat}
  =\D{Y}\tx c^*
\,,\]
\DE{}\EN{and of course $m=m^*$.}
\end{prf}

\DE{}
\EN{This transformation rule can look quite complex if written in single terms.
On the other hand, the general description can easily be remembered and is motivated by the following representation}
\begin{equ}{vv}
  F_{k_1\cdots k_M}
  =\rho{\tx v}_{k_1}\cdots{\tx v}_{k_M}+\tx\Pi_{k_1\cdots k_M}
\DE{}\EN{\,.}\end{equ}
\DE{}
\EN{We remark that the transformation rule \EQU{trans.o} works
also for arbitrary transformations $Y$, hence its a rule which
we will postulate also in the relativistic case, see \EQU{higher.*.rule-T}.}

\medskip
\DE{}
\EN{Furthermore, the transformation rule \EQU{trans.o} is important
in connection with the general rule \EQU{continuum.*.o}.}
\DE{}
\EN{We define the physical properties of the quantities in \EQU{moments-dgl}
by saying that in the weak formulation}
\begin{Equ}{moments-weak}
  \int_{\RR^4}\Big(
  \sum_{k=0}^3\de_{y_k}\zeta_{i_1\cdots i_N}\cdot F_{i_1\cdots i_Nk}
  + \zeta_{i_1\cdots i_N}R_{i_1\cdots i_N}
  \Big)\d\leb{4}=0
\end{Equ}
\DE{}
\EN{the test functions $\zeta_{i_1\cdots i_N}\in C_0^\infty(\RR^4)$
for $i_1\til i_N\in\{0,1\til 3\}$
satisfy the following transformation rule}
\[
  \zeta^*_{\bar i_1\cdots\bar i_N}
  =\sum_{i_1\til i_N=0}^3
  Y_{i_1\p\bar i_1}\cdots Y_{i_N\p\bar i_N}
  \zeta_{i_1\cdots i_N}\circ Y
\]
\DE{}
\EN{for $\bar i_1\til\bar i_N\in\{0,1\til 3\}$}
\DE{}
\EN{This means that $\zeta^*=\up{Z}T\zeta\circ Y$ where}
\[
  Z_{(i_1\til i_N)(\bar i_1\til\bar i_N)}
  =
  Y_{i_1\p\bar i_1}\cdots Y_{i_N\p\bar i_N}
\,.\]
\DE{}
\EN{By \REF{continuum.thm.} this is satisfied if}
\begin{Equ}{F}
  F_{i_1\cdots i_Nk}\circ Y
  =\sum_{\bar i_1\til\bar i_N,\bar k=0}^n
  Z_{(i_1\til i_N)(\bar i_1\til\bar i_N)}
  Y_{k\p\bar k}
  F^*_{\bar i_1\cdots\bar i_N\bar k}
\end{Equ}
\DE{}\EN{and}
 \begin{Equ}{g}
  R_{i_1\cdots i_N}\circ Y
  =\sum_{\bar i_1\til\bar i_N,\bar k=0}^n
  (Z_{(i_1\til i_N)(\bar i_1\til\bar i_N)})_{\p\bar k}
  F^*_{\bar i_1\cdots\bar i_N\bar k}
\\\hfill
  +\sum_{\bar i_1\til\bar i_N}^n
  Z_{(i_1\til i_N)(\bar i_1\til\bar i_N)}
  R^*_{\bar i_1\cdots\bar i_N}
\,.\end{Equ}
\DE{}
\EN{Equation \EQU{F} is equivalent to \EQU{trans.o} for $M=N+1$
which was proved in \REF{trans.}. The proof of \EQU{g} you will find
in \CIT{Alt-Kontinuum}{Chap V}.}
\DE{}
\EN{This presentation will serve us in section~\REF{higher..}
as guide.}

{}
\newpage
{{}
\sect{time}{\DE{}\EN{Time and space}}

\DE{}\EN{The points in spacetime are denoted by $y=(y_0,y_1,y_2,y_3)\in\RR^4$ and
in spacetime a symmetric $4\times4$-matrix $y\mapsto\GG(y)$ is given.
This matrix describes the hyperbolic geometry, that is,
\begin{equ}{lambda}\begin{arr}c
  \text{$\GG$ has a negative eigenvalue: $\lambda_0<0$,}
  \\\text{the remaining eigenvalues of $\GG$ are positive:
  $\lambda_i>0$ for $i\geq1$.}
\end{arr}\end{equ}
If follows from the theory of symmetric matrices:
}
\begin{stmt}{orth}{Theorem}
\DE{}
\EN{Let $\lambda_k$, $k\geq0$, be the eigenvalues of $\GG$ as in \EQU{*.lambda},
then there exists a orthonormal basis $\{e^\orth_0\til e^\orth_n\}$ with}
\[
  \GG=\sum_{k\geq0}\lambda_ke^\orth_k\up{(e^\orth_k)}T
  =\sum_{k\geq0}\lambda_ke^\orth_k\otimes e^\orth_k
\,.\]
\DE{}
\EN{The eigenfunctions $e^\orth_k$ depend on $y$.}
\end{stmt}
\begin{prf}{}\DE{}\EN{It is} $\GG e^\orth_k=\lambda_ke^\orth_k$
\DE{}\EN{and} $e^\orth_k\dd e^\orth_l=\kronecker{kl}$.
\end{prf}

\DE{}\EN{We assume that to every point $y\in\RR^4$ there
exists a (dual) direction $\ewelt(y)\neq0$, which we call ``direction of time'',
such that}\footnote{\DE{}\EN{We denote by ``$\dd$'' the inner product in spacetime.
}}
\begin{Equ}{zeit}
  \ewelt\dd\GG\ewelt=-\frac1{\cc^2}
\,.\end{Equ}
\DE{}
\EN{The ``speed of light'' $\cc>0$ occurs here in \EQU{zeit}
for the first time, and is the same for all observers. This follows from
the fact that \EQU{zeit} is objective as shown in \REF{change.cons.}.
The ``space'' which is orthogonal to $\ewelt$ we denote by}
\begin{equ}{welt}
  \welt(y):=\up{\{\ewelt(y)\}}\perp
  =\set{z\in\RR^4}{\ewelt(y)\dd z=0}
\,.\end{equ}
\DE{}
\EN{This space $\welt(y)\subset\RR^4$ is $3$-dimensional,
orthogonal to $\{\ewelt(y)\}$,
and $\ewelt(y)$ is a point in $\RR^4$.}

\medskip
\DE{}\EN{It turns out that $\ewelt$ has to depend on $y$
(this is related to the theorem on second derivatives of
observer transformations),
hence $\ewelt$ in general is not constant.
Now $(\ewelt,\welt)$ describe the coordinates of the world
for a single observer. In the next section we see what this description
says for another observer. In the special case $\GG=\GG_\cc$ we have
the situation of Lorentz observers, and this case serves as an example.
}\DE{}
\EN{It is here and hereinafter always $n=3$.}

\begin{stmt}{standard}{\DE{}\EN{Standard example}}
It is $\GG=\GG_\cc$, where
\begin{Equ}{o}
  \GG_\cc=-\frac1{\cc^2}\ee_0\up{\ee_0}T+\sum_{i\geq1}\ee_i\up{\ee_i}T
  =\begin{mat}-\frac1{\cc^2}&0\\0&\Id\end{mat}
\,.\end{Equ}
We then have as an example
\[
  \ewelt:=\ee_0 \quad\text{and}\quad \welt:=\fcn{span}\{\ee_1\til\ee_n\}
\,.\]
Here $\{\ee_0\til\ee_n\}$ is the standard orthonormal basis of $\RR\times\RR^n$.
The classical limit one obtains for $\cc\to\infty$.
\end{stmt}

\DE{}
\EN{About $\ewelt$ the following statement is true.}

\begin{stmt}{dual}{Lemma}
\DE{}\EN{Let} $e_0':=\ewelt$ \DE{}\EN{and let} $\{e_1\til e_n\}$
\DE{}\EN{be a basis of} $\welt:=\up{\{\ewelt(y)\}}\perp$,
\DE{}\EN{then there exists for} $(\nu_1\til\nu_n)$
\DE{}\EN{one and only one}
$e_0$ \DE{}\EN{and} $e_i'$ \DE{}\EN{for} $i=1\til n$
\DE{}\EN{with} $e_i'+\nu_ie_0'\in\welt$, \DE{}\EN{so that}
$\{e_0\til e_n\}$ \DE{}\EN{and} $\{e_0'\til e_n'\}$
\DE{}
\EN{are dual basis, i.e.}
\[
  e_k'\dd e_l=\kronecker{k,l} \tm{\DE{}\EN{for}} k,l=0\til n
\,.\]
\begin{rem}{\DE{}\EN{Hint}}
\DE{}
\EN{The free parameters $(\nu_1\til\nu_n)$ correspond to a
``velocity'' of an observer transformation, see \REF{unique.}.}
\end{rem}
\end{stmt}
\begin{prf}{}
\DE{}
\EN{Using \EQU{*.welt} it follows that $\{e_0',e_1\til e_n\}$
is a basis of spacetime. Let with certain coefficients}
\[
  e_0:=\mu e_0' + \sum_{j=1}^n\nu_je_j
\]
\DE{}
\EN{be a vector. The property $e_0'\dd e_0=1$ of a dual matrix gives}
\[
  1=e_0'\dd e_0=\mu |e_0'|^2
\quad\text{\DE{}\EN{hence}}\quad
  \mu=\frac1{|e_0'|^2}
\,.\]
\DE{}
\EN{Next let with certain coefficients}
\[
  e_i':=b_i e_0' + \sum_{k=1}^na_{ik}e_k \tm{\DE{}\EN{for}} i=1\til n
.\]
\DE{}
\EN{Then we obtain from the property of a dual basis for $i,j=1\til n$}
\[\begin{arr}c
  \kronecker{i,j}=e_i'\dd e_j=\sum_{k=1}^na_{ik}e_k\dd e_j
  =(AE)_{ij}
\\
  \text{\DE{}\EN{where}~} A:=\seq{a_{ik}}{i,k=1\til n}
  \text{~\DE{}\EN{and}~} E:=\seq{e_k\dd e_j}{k,j=1\til n}
,\end{arr}\]
\DE{}
\EN{hence $\Id=AE$, that is $A=\up{E}{-1}$. And for $i=1\til n$}
\[\begin{arr}l
  0=e_i'\dd e_0=\Big(b_i e_0' + \sum_{k=1}^na_{ik}e_k\Big)
  \dd\Big(\mu e_0' + \sum_{j=1}^n\nu_je_j\Big)
\\
  =b_i\mu|e_0'|^2+\sum_{j,k=1}^na_{ik}\nu_je_k\dd e_j
  =b_i+\sum_{j=1}^n\kronecker{i,j}\nu_j=b_i+\nu_i
\,,\end{arr}\]
\DE{}\EN{hence} $b_i=-\nu_i$.
\end{prf}

\DE{}
\EN{Two every basis there exists one and only one dual basis, this is a
simple consequence of Functional Analysis.}
\DE{}
\EN{Now let $\{e_0\til e_n\}$ and $\{e_0'\til e_n'\}$ be any dual basis.
If one defines $\ewelt:=e_0'$ and $W:=\fcn{span}\{e_1\til e_n\}$
then it follows \EQU{welt}, however, the property \EQU{zeit}
is still to be satisfied.
For given $\ewelt$ with \EQU{zeit} one can define a dual basis
also as follows.}

\begin{stmt}{thm}{Theorem}
\DE{}
\EN{Let $\ewelt$ with \EQU{*.zeit} and let $\welt:=\up{\{\ewelt\}}\perp$.
We assume that}
\[
  (z,w)\in\welt\times\welt \quad\longmapsto\quad z\dd\up{\GG}{-1}w \in\RR
\]
\DE{}
\EN{is a scalar product, i.e.~$w\dd\up{\GG}{-1}w>0$ for
$w\in\welt\setminus\{0\}$.}
\begin{enum}
\num{e0}
\DE{}
\EN{Choose $e_0':=\ewelt$ and define $e_0:=-\cc^2\GG e_0'$.}
\num{W}
\DE{}
\EN{Choose a basis $\{e_1\til e_n\}$ of $\welt$
which is $\up\GG{-1}$-orthogonal, i.e.
\[e_j\dd\up\GG{-1}e_i=\kronecker{i,j}\tm{f\"ur}i,j=1\til n.\]
\vspace*{-7mm}}
\num{ei}
\DE{}
\EN{Define $e_i':=\up\GG{-1}e_i$ for $i=1\til n$.}
\end{enum}
\DE{}
\EN{Then $\{e_0\til e_n\}$ and $\{e_0'\til e_n'\}$ are dual basis
and}
\begin{Equ}{o}  \GG=-\frac1{\cc^2} e_0\up{e_0}T+\sum_{i=1}^ne_i\up{e_i}T
\,.\end{Equ}\end{stmt}

\begin{picture}(100,274)
\put(40,0){{\includegraphics[height=10cm]{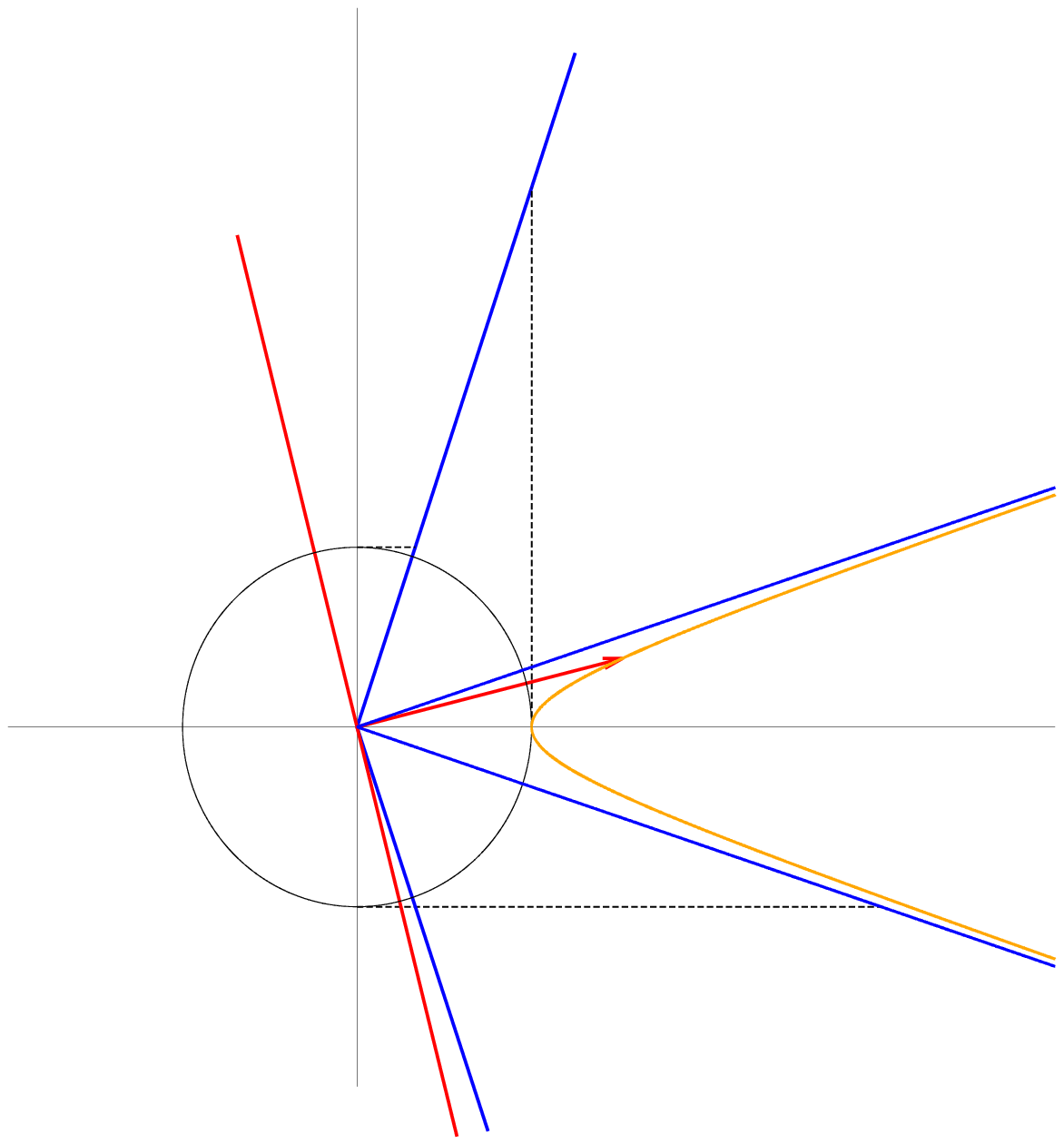}}}
{\color{red}
\put(80,225){$\welt=\up{\{\ewelt\}}\perp$}
\put(200,105){$\ewelt$}
}{\color{orange4}
\put(240,115){$\set{z'}{z'\dd\GG z'=-\frac1{\cc^2}}$}
\color{black}\put(230,116){$\nwarrow$}
}{\color{blue}
\put(200,155){$\set{z'}{z'\dd\GG z'=0}$}
\color{black}\put(254,145){$\downarrow$}
}{\color{blue}
\put(201,255){$\set{z}{z\dd\up\GG{-1}z=0}$}
\color{black}\put(191,255){$\leftarrow$}
}{\color{black}
\put(144,151){$\displaystyle\frac1\cc$}
\put(183,200){$\cc$}
\put(193,47){$\cc$}
}\end{picture}

\begin{prf}{\DE{}\EN{of the duality}}
\DE{}\EN{It is} $e_0'\dd e_0=-\cc^2e_0'\dd\GG e_0'=1$.
\DE{}\EN{Since} $e_i\in\welt=\up{\{e_0'\}}\perp$
\DE{}\EN{it is} $e_0'\dd e_i$ =0
\DE{}\EN{for} $i=1\til n$. \DE{}\EN{And}
$
  e_i'\dd e_j=e_j\dd\up\GG{-1}e_i=\kronecker{i,j}
$
\DE{}\EN{for} $i,j=1\til n$ \DE{}\EN{by construction.}
\DE{}
\EN{Since $\GG^{-1}$ is symmetric it follows}
\[
  e_i'\dd e_0=e_0\dd e_i'
  =e_0\dd(\up\GG{-1}e_i)=(\up\GG{-1}e_0)\dd e_i=-\cc^2e_0'\dd e_i=0
\]
\DE{}\EN{for} $i=1\til n$.
\end{prf}
\begin{prf}{\DE{}\EN{of the representation of} $\GG$}
Define
\[
  \tilde\GG:=-\frac1{\cc^2}e_0\up{e_0}T+\sum_{i=1}^ne_i\up{e_i}T
\,.\]
The dual basis implies
$\tilde\GG e_0'=-\frac1{\cc^2}e_0=\GG e_0'$ and $\tilde\GG e_i'=e_i=\GG e_i'$ for $i=1\til n$. Hence $\tilde\GG=\GG$.
\end{prf}

\DE{}\EN{If $e_0'$ as in \REF{thm.e0} and $\{e_1\til e_n\}$ as in \REF{thm.W}
are chosen then we can represent $e_0$ as in the proof of \REF{dual.}
\[
  e_0=\mu e_0' + \sum_{j=1}^n\nu_je_j \,,\quad \mu=\frac1{|e_0'|^2}
\,.\]
If we define $e_i':=\up\GG{-1}e_i$ for $i=1\til n$ as in \REF{thm.ei}
then
}\[\begin{arr}l
  0=e_j'\dd e_0=e_0\dd\up\GG{-1}e_j
  =\frac1{|e_0'|^2}e_0'\dd\up\GG{-1}e_j+\sum_i\nu_ie_i\dd\up\GG{-1}e_j
\\
  =\frac1{|e_0'|^2}e_0'\dd\up\GG{-1}e_j+\nu_j
\quad\tm{\DE{}\EN{hence}}\quad
  \nu_j=-\frac{e_0'\dd\up\GG{-1}e_j}{|e_0'|^2}=-\frac{e_0'\dd e_j'}{|e_0'|^2}
.\end{arr}\]
\DE{}\EN{Consequently, the freedom in the choice of $(\nu_1\til\nu_n)$
is given by the choice of the basis $\{e_1\til e_n\}$ in \REF{thm.e0}.
In \REF{thm.} the following properties of the matrix $\GG$ are addressed,
where \REF{prop.W} is an essential property, whereas \REF{prop.G} is
important for practical use.
}

\begin{stmt}{prop}{Properties}
\DE{}
\EN{For a symmetric matrix $\GG$ consider the following properties:}
\begin{enum}
\num{G}
\DE{}
\EN{It is $\{e_0\til e_n\}$ a basis and}
\[
  \GG=-\frac1{\cc^2} e_0\up{e_0}T+\sum_{i=1}^ne_i\up{e_i}T
\,.\]
\vspace*{-7mm}
\num{eigen}
\DE{}
\EN{It is $\{e_0\til e_n\}$ a basis and
$\{e_0'\til e_n'\}$ the corresponding dual basis and}
\[
  \GG e_0'=-\frac1{\cc^2} e_0
\,,\quad
  \GG e_i'=e_i\tm{f\"ur}i\geq1
\,.\]
\vspace*{-7mm}
\num{W}
\DE{}
\EN{It is $\welt\subset\RR^4$ a subspace of codimension 1 and}
\[\text{
  \DE{}
  \EN{$(z,w)\mapsto z\dd\up{\GG}{-1}w$ on $\welt$ a scalar produkt.} 
}\]
\end{enum}
\end{stmt}
\DE{}
\EN{The connection between these properties is the content of the following
lemmata.}
\DE{}
\EN{The property \REF{prop.G} implies immediately,
if $\{e_0'\til e_n'\}$ is the dual basis,}
\begin{Equ}{-}
  \up\GG{-1}=-\cc^2 e_0'\up{e_0'}T+\sum_{i=1}^ne_i'\up{e_i'}T
\,.\end{Equ}
\DE{}
\EN{Let $\ewelt=e'_0$. The property \REF{prop.eigen} implies immediately
\EQU{zeit}, since $\GG e'_0=-\frac1{\cc^2}e_0$ implies}
\[
  e'_0\dd\GG e'_0=-\frac1{\cc^2}e'_0\dd e_0=-\frac1{\cc^2}
\,.\]

\begin{stmt}{lem1}{Lemma}
\DE{}
\EN{\REF{prop.G} and \REF{prop.eigen} are equivalent.}
\end{stmt}
\DE{}
\EN{That \REF{prop.eigen} implies \REF{prop.G} has been proved in \REF{thm.}.}
\begin{prf}{\REF{prop.G}$\Rightarrow$\REF{prop.eigen}}
\DE{}
\EN{With $\lambda_0=-\frac1{\cc^2}$ and $\lambda_i=1$ for $i\geq1$ it implies,
if we define $e_k':=\lambda_k\up\GG{-1}e_k$,}
\[
  \lambda_le_l=\GG e_l'=\big(\sum_{k\geq0}\lambda_ke_k\up{e_k}T\big)e_l'
  =\sum_{k\geq0}\lambda_k(e_k\dd e_l')e_k
\,.\]
\DE{}\EN{Hence it holds for all} $k,l\geq0$
\[
  \lambda_ke_k\dd e_l'=\lambda_l\kronecker{k,l}=\lambda_k\kronecker{k,l}
\,.\]
\DE{}
\EN{Since all $\lambda_k\neq0$ we conclude}
$
  e_k\dd e_l'=\kronecker{k,l}
$.
\DE{}
\EN{This says that $\set{e_l'}{l\geq0}$ is the dual basis.}
\end{prf}

\begin{stmt}{lem2}{Lemma}
\DE{}
\EN{Let \REF{prop.G} be true where $\{e'_0\til e'_n\}$ is the dual basis
of $\{e_0\til e_n\}$. If $\ewelt=e_0'$ and $\welt=\up{\{e_0'\}}\perp$
then the property \REF{prop.W} is true.}
\end{stmt}

\begin{prf}{}
\DE{}
\EN{For $z'$ and $z$ we have the representation}
\[\begin{arr}l
  z'=\sum_{k\geq0}z'_ke'_k\,,\quad z'_k=z'\dd e_k
\,,\\
  z=\sum_{k\geq0}z_ke_k\,,\quad z_k=z\dd e'_k
\end{arr}\]
\DE{}\EN{and we have $\ewelt=e'_0$ and}
$\welt=\up{\{e'_0\}}\perp=\fcn{span}\{e_1\til e_n\}$.
\DE{}\EN{From \REF{prop.G} it follows}
\[\begin{arr}l
  \GG z'=-\frac1{\cc^2}z'_0e_0+\sum_{i\geq1}z'_ie_i
\,,\\
  z'\dd\GG z'=-\frac1{\cc^2}|z'_0|^2+\sum_{i\geq1}|z'_i|^2
\,.\end{arr}\]
\DE{}
\EN{If $z'=\up\GG{-1}z$, that means $\GG z'=z$, then}
\[
  -\frac1{\cc^2}z'_0e_0+\sum_{i\geq1}z'_ie_i=z=\sum_{k\geq0}z_ke_k
\,,\]
\DE{}\EN{and therefore}
\[
  z_0=-\frac1{\cc^2}z'_0\,,\quad z_i=z'_i\tm{\DE{}\EN{for}}i\geq1
.\]
\DE{}
\EN{Now let $z\in\welt$, that is $z_0=z\dd e'_0=0$, and then also $z'_0=0$.
This implies}
\[
  z\dd\up\GG{-1}z=(\GG z')\dd z'=z'\dd\GG z'
  =\sum_{i\geq1}|z'_i|^2\geq0
\EN{\,.}\]
\DE{}
\EN{And this is strict positive if $z\neq0$ which is equivalent to $z'\neq0$.}
\end{prf}

\DE{}
\EN{The following lemma shows that for all occurring matrices $\GG$ the
property \REF{prop.W} is satisfied, and therefore also the construction
in \REF{thm.}.}

\begin{stmt}{lem3}{Lemma}
\DE{}
\EN{Let $\GG$ be a matrix as in \EQU{*.lambda} and
$\ewelt$ with \EQU{*.zeit} and $\welt=\up{\{\ewelt\}}\perp$.}
\DE{}
\EN{Then the property \REF{prop.W} is true.}
\end{stmt}

\begin{prf}{}
\DE{}
\EN{Let $\lambda_0<0$ and $\lambda_i>0$ for $i\geq1$.}
For vectors $z'$ one has the identity
\[
  \GG z'=\sum_{k\geq0}\lambda_kz'_ke^\perp_k
\,,\quad
  z'_k:=z'\dd e^\perp_k
\]
and therefore
\begin{Equ}{1}
  z'\dd\GG z'=\sum_{k\geq0}\lambda_k|z'_k|^2
\,.\end{Equ}
And it follows that for vectors $z$ one has the identity
\[
  \up\GG{-1}z=\sum_{k\geq0}\lambda_k^{-1}z_k e^\perp_k
\,,\quad
  z_k:=z\dd e^\perp_k
\]
hence
\begin{Equ}{2}
  z\dd\up\GG{-1}z=\sum_{k\geq0}\frac{|z_k|^2}{\lambda_k}
\,.\end{Equ}
Since $\ewelt\dd\GG\ewelt=-\frac1{\cc^2}$ by \EQU{*.zeit} and hence
\[
  \bar\ewelt\dd\GG\bar\ewelt=\lambda_0=-|\lambda_0|
  \,,\quad \bar\ewelt:=\cc\sqrt{|\lambda_0|}\,\ewelt
\,.\]
The first identity \EQU{1} gives (with $z'=\bar\ewelt$)
\[
  -|\lambda_0|=-|\lambda_0|\cdot|\bar\ewelt_0|^2
  +\sum_{i\geq1}\lambda_i|\bar\ewelt_i|^2
\,,\quad
  \bar\ewelt_k:=\bar\ewelt\dd e^\perp_k
\,,\]
that is
\begin{Equ}{e}
  |\bar\ewelt_0|^2=1+\frac1{|\lambda_0|}\sum_{i\geq1}\lambda_i|\bar\ewelt_i|^2
\,.\end{Equ}
Therefore $\bar\ewelt_0\neq0$ and for
$
  z\in\welt=\up{\{\ewelt\}}\perp\setminus\{0\}
  =\up{\{\bar\ewelt\}}\perp\setminus\{0\}
$, that is,
\[
  0=z\dd\bar\ewelt=z_0\bar\ewelt_0+\sum_{i\geq1}z_i\bar\ewelt_i
\quad\text{hence}\quad
  -z_0=\sum_{i\geq1}z_i\frac{\bar\ewelt_i}{\bar\ewelt_0}
\,,\]
we get
\[
  |z_0|^2\leq\Big(\sum_{i\geq1}|z_i|\frac{|\bar\ewelt_i|}{|\bar\ewelt_0|}\Big)^2
  \leq\sum_{i\geq1}\frac{|z_i|^2}{\lambda_i}\cdot
  \sum_{i\geq1}\lambda_i\frac{|\bar\ewelt_i|^2}{|\bar\ewelt_0|^2}
\,,\]
where from \EQU{e}
\[
  {\sum_{i\geq1}\lambda_i\frac{|\bar\ewelt_i|^2}{|\bar\ewelt_0|^2}}<|\lambda_0|
\,.\]
Then from the second identity \EQU{2}
\[
  z\dd\up\GG{-1}z
  =-\frac{|z_0|^2}{|\lambda_0|}+\sum_{i\geq1}\frac{|z_i|^2}{\lambda_i}
  \geq\sum_{i\geq1}\frac{|z_i|^2}{\lambda_i}\cdot
  \Big(1-\frac1{|\lambda_0|}
  \sum_{i\geq1}\lambda_i\frac{|\bar\ewelt_i|^2}{|\bar\ewelt_0|^2}\Big)
  >0
\]
which had to be shown.
\end{prf}

\DE{}
\EN{That for given $\ewelt$ the dual basis has free parameters
$(\nu_1\til\nu_n)$ we have seen in the proof of \REF{dual.}.}
\DE{}
\EN{Therefore the representation \REF{prop.G} for $\GG$ is not unique.
We see this at the standard matrix $\GG_\cc$ in \EQU{standard.o}:}

\begin{stmt}{unique}{\DE{}\EN{Theorem}} \DE{}\EN{Let} $\GG=\GG_\cc$ \DE{}\EN{and let} $\{e_0\til e_n\}$
\DE{}
\EN{a (on $y$ depending) basis. Then} 
\[
  \GG_\cc=-\frac1{\cc^2} e_0\up{e_0}T+\sum_{i=1}^ne_i\up{e_i}T
\]
\DE{}
\EN{if and only if modulo the sign of each basis vector
there exists a Lorentz matrix ${\bf L}_c(V,Q)$
(here $V$ and $Q$ depend on $y$) with}
\[\begin{arr}l
 e_0=\begin{mat}\gamma\\\gamma V\end{mat}
,\quad
  e_0'=-\frac1{\cc^2}\up{\GG_\cc}{-1}e_0
  =\begin{mat}\gamma\\-\frac\gamma{\cc^2}V\end{mat}
,\quad
  \tm{\DE{}\EN{and for}}i\geq1
:\\
  e_i=\begin{mat}\frac\gamma{\cc^2}V\dd Q\ee_i  \\Q\ee_i+\frac{\gamma^2V\dd Q\ee_i}{\cc^2(\gamma+1)} V\end{mat}
,\quad
  e_i'=\up{\GG_\cc}{-1}e_i
  =\begin{mat}-\gamma V\dd Q\ee_i
  \\Q\ee_i+\frac{\gamma^2V\dd Q\ee_i}{\cc^2(\gamma+1)} V\end{mat}
\,.\end{arr}\]
\DE{}
\EN{And of course $e'_k\dd e_l=\kronecker{kl}$ for $k,l\geq0$.}
\end{stmt}

\begin{prf}{}
\DE{}
\EN{Let us write the vectors $e_k$ in components
$e_k=(M_{1k}\til M_{nk})$, so that the matrix $M=\seq{M_{ij}}{ij}$
satisfies with the canonical basis vectors $\ee_i$
the equation $M_{ik}=e_k\dd\ee_i$.
Then $e_k=M\ee_k$ and with 
$\lambda_0:=-\frac1{\cc2}$ and $\lambda_i:=1$ for $i\geq1$}
\[\begin{arr}l
  \GG_\cc=\sum_{k\geq0}\lambda_ke_k\up{e_k}T
  =\sum_{k\geq0}\lambda_kM\ee_k\up{(M\ee_k)}T
\\
  =M\Big(\sum_{k\geq0}\lambda_k\ee_k\up{\ee_k}T\Big)\up{M}T=M\GG_\cc\up{M}T
\,.\end{arr}\]
This says that the Matrix $M$ keeps $\GG_\cc$ unchanged, and
therefore implies (see \REF{invariance.thm.}
\DE{}\EN{under the assumption $M_{00}\geq0$ and $\fcn{det}M\geq0$})
that $M$ is a Lorentz-Matrix. \DE{}
\EN{The rows of this matrix are $e_k=M\ee_k$ for $k\geq0$.}
\end{prf}

\DE{}
\EN{As generalization we give in \REF{change.general.}
representations of $\GG$ for different $\ewelt$-vectors.}

\begin{stmt}{klass}{\DE{}\EN{Remark (Classical physics)}}
\DE{}
\EN{In the classical limit $\cc\to\infty$ the basis
$\{e_0\til e_n\}$ and $\{e_0'\til e_n'\}$ in \REF{unique.} converge to}
\[\begin{arr}l
 e_0=\begin{mat}1\\V\end{mat}
,\quad
  e_0'=\begin{mat}1\\0\end{mat}
,\quad
  \tm{\DE{}\EN{and for}}i\geq1
:\\
  e_i=\begin{mat}0\\Q\ee_i\end{mat}
,\quad
  e_i'=\begin{mat}-V\dd Q\ee_i\\Q\ee_i\end{mat}
\,,\end{arr}\]
and the matrix is
\[
  \GG_\infty=\sum_{i\geq1}e_i\up{e_i}T
\,.\]
\end{stmt}

{}

{}

\sect{change}{\DE{}\EN{Change of observer}}

\DE{}
\EN{Since the vector $\ewelt$ will occur in the differential equations
we have to guarantee the rule by which this quantity will change
between observers. If $\ewelt$ is this vector for one observer
and $\ewelt^*$ is this vector for another observer, the transformation rule is}
\begin{equ}{trans}
  \ewelt^*=\up{(\tx\D Y)}T\,\ewelt\circ Y
\,,\end{equ}
\DE{}
\EN{where $Y$ is the observer transformation.
This means that $\ewelt$ is a covariant vector.}

\begin{stmt}{cons}{\DE{}\EN{Consistence}}
\DE{}
\EN{The transformations rule \EQU{*.trans} is consistent with
the assumption \EQU{time.*.zeit} and implies}
\[\welt\circ Y=\D{Y}\,\welt^*\,.\]
\end{stmt}
\begin{prf}{}\DE{}\EN{It is}
\[\begin{arr}l
  (\ewelt\dd\GG\ewelt)\circ Y=(\ewelt\circ Y)\dd
  (\tx\D{Y}\GG^*\up{\tx\D{Y}}T\ewelt\circ Y)
\\
  =(\up{\tx\D{Y}}T\ewelt\circ Y)\dd\GG^*\up{\tx\D{Y}}T\ewelt\circ Y
  =\ewelt^*\dd\GG^*\ewelt^*
\,,\end{arr}\]
\DE{}
\EN{hence the condition $\ewelt\dd\GG\ewelt=-\frac1{\cc^2}$ is objective.}
\DE{}
\EN{And for $w^*$ with $w:=\D{Y}w^*$ it holds}
\[
  (\ewelt\dd w)\circ Y=(\ewelt\circ Y)\dd\D{Y}w^*
  =(\up{\D{Y}}T\ewelt\circ Y)\dd w^*=\ewelt^*\dd w^*
\,,\]
\DE{}
\EN{hence the condition $\ewelt\dd w=0$, which defines $\welt$, is objective.}
\end{prf}

\DE{}
\EN{The transformation rule for $\ewelt$ can be generalized to the
basis elements, where again $n=3$.} 

\begin{stmt}{basis-trans}{\DE{}\EN{Consistence with basis} (Definition)}
\DE{} \EN{Let $\{e_0\til e_n\}$ be a basis with dual basis
$\{e_0'\til e_n'\}$.}
\DE{}
\EN{Then the transformation rules for $k\geq0$ are}
\begin{equ}{o}\begin{arr}c
  e_k\circ Y=\tx\D{Y}e_k^*
\,,\\
  e_k'^*=\up{(\tx\D{Y})}Te_k'\circ Y
\,.\end{arr}\end{equ}
\DE{}
\EN{Setting $\ewelt=e'_0$ this is in accordance with \EQU{*.trans}.}
\DE{}
\EN{The definition shows: The basis elements $e_k$ are contravariant vectors and
the dual elements $e'_k$ covariant vectors.}
\\\tit{\DE{}\EN{Assertion}}
\DE{}
\EN{The transformation rules are compatible with the definition of
a dual basis.}
\end{stmt}
\begin{prf}{\DE{}\EN{of the assertion}}
\DE{}\EN{If the transformation rules are true for $\{e_0\til e_n\}$
and if the dual basis $\{e_0'\til e_n'\}$
for one observer is given, then it holds}
\[
  \kronecker{k,l}=(e_k'\dd e_l)\circ Y=(e_k'\circ Y)\dd(\tx\D{Y}e_l^*)
  =(\up{\tx\D{Y}}Te_k'\circ Y)\dd e_l^*
\DE{}\EN{\,.}\]
\DE{}\EN{Hence $e_k'^*:=\up{\tx\D{Y}}Te_k'\circ Y$ is the dual basis
of $\{e_0^*\til e_n^*\}$.}
\DE{}\EN{The other way around, if the transformation rules are true for
$\{e_0'\til e_n'\}$ and if the dual basis $\{e_0\til e_n\}$
for one observer is given, then it holds}
\[
  \kronecker{k,l}=(e_k'\dd e_l)\circ Y
  =(\up{(\tx\D{Y})}{-T}e_k'^*)\dd e_l\circ Y
  =e_k'^*\dd(\up{(\tx\D{Y})}{-1}e_l\circ Y)
\DE{}\EN{\,.}\]
\DE{}\EN{Hence $e_l^*:=\up{\tx\D{Y}}{-1}e_l\circ Y$ is the dual basis
of $\{e_0'^*\til e_n'^*\}$.}
\end{prf}

\DE{}\EN{It follows from the transformation rules \EQU{basis-trans.o}
that $e_k'\dd e_l$ are objective scalars.}

\begin{stmt}{G*}{\DE{}\EN{Lemma}}
\DE{}\EN{Let for the other observer}
\[
  \GG^*:=-\frac1{\cc^2} e_0^*\up{e_0^*}T+\sum_{i=1}^ne_i^*\up{e_i^*}T
\]
\DE{}\EN{and let $y^*\mapsto y=Y(y^*)$ be the observer transformation.}
\DE{}\EN{Then it holds for the local observer}
\[
  \GG=-\frac1{\cc^2} e_0\up{e_0}T+\sum_{i=1}^ne_i\up{e_i}T
\,,\]
\DE{}\EN{if the basis is transformed according to \EQU{basis-trans.o}.}
\end{stmt}
\begin{prf}{}
\DE{}\EN{Let $\lambda_0=-\frac1{\cc^2}$ and $\lambda_i=1$ for $i\geq1$. Then}
\[\begin{arr}l
  \GG\circ Y=\tx\D{Y}\GG^*\up{(\tx\D{Y})}T
  =\sum_{k\geq0}\lambda_k\tx\D{Y}e_k^*\up{e_k^*}T\up{(\tx\D{Y})}T
\\
  =\sum_{k\geq0}\lambda_k(\tx\D{Y}e_k^*)\up{(\tx\D{Y}e_k^*)}T
  =\Big(\sum_{k\geq0}\lambda_ke_k\up{e_k}T\Big)\circ Y
\,,\end{arr}\]
\DE{}\EN{since $e_k\circ Y=\tx\D{Y}e_k^*$ by \EQU{basis-trans.o}.}
\end{prf}

\DE{}
\EN{We remark that all observer transformations are allowed
which convert the hyperbolic geometry again in a geometry of the same type.
Since they are connected with the standard situation it follows
from the group property that the general matrix $\GG$
can be expressed by $\GG_\cc$ and general transformations. 
Therefore we have the following theorem, where we assume that an observer
has the matrix $\GG$ and there exists an observer transformation 
$Y$ to a standard observer.}

\begin{stmt}{general}{Theorem}
\DE{}\EN{Let $\{e_0\til e_n\}$ be an on $y$ depending basis
and $\{e_0'\til e_n'\}$ the corresponding dual basis
with $e_0'\dd\GG e_0'=-\frac1{\cc^2}$.
Then there exists modulo the sign of the basis elements
a Lorentz matrix ${\bf L}_c(V,Q)$,
which depends on $y$, with}
\[
  e_k\circ Y=\D{Y}{\bf L}_c(V,Q)\ee_k
\,,\quad
  \up{\big(\D{Y}{\bf L}_c(V,Q)\big)}Te_k'\circ Y=\ee_k
\,.\]
\DE{}\EN{Here $\ee_k$ are as in \REF{time.standard.}.}
\end{stmt}
\begin{prf}{}
\DE{}\EN{There exists an observer transformation $y=Y(y^*)$ such that
$y^*$ are the coordinates of the standard observer.
Let $e_k^*$ and $e_k'^*$ be the corresponding basis vectors, that is,}
\[
  e_k\circ Y=\tx\D{Y}e_k^*
\,,\quad
  e_k'^*=\up{(\tx\D{Y})}Te_k'\circ Y
\,.\]
\DE{}\EN{Then $\{e_0^*\til e_n^*\}$ and $\{e_0'^*\til e_n'^*\}$
are dual basis of the standard observer which satisfy
$e_0'^*\dd\GG_\cc e_0'^*=-\frac1{\cc^2}$, hence by \REF{time.unique.} 
modulo the sign of the elements there exists
a Lorentz matrix ${\bf L}_c(V,Q)$ with}
\[
  e_k^*={\bf L}_c(V,Q)\ee_k
\,,\quad
  \up{{\bf L}_c(V,Q)}Te_k'^*=\ee_k
\,.\]
This implies the assertion.
\end{prf}

{}}
\newpage
{{}
\sect{scalar}{\DE{}\EN{Scalar conservation laws}}
We introduce here the simplest example of conservations laws,
which is a single equation
\begin{equ}{strong}
  \tx{\div}\,q=\rate
\end{equ}
in a domain $\UU\subset\RR^4$, where $q=\set{q_k}{k=0\til 3}$
is a 4-vector and $\rate$ an objective scalar.
The weak version of this law reads \begin{Equ}{weak}
  \int_\UU\big( 
    \tx{\grad}\eta\dd q + \eta\rate
  \big)\d\leb4=0
\end{Equ}
for test functions $\eta\in\C_0^\infty(\UU;\RR)$.
The equation \EQU{*.strong} is called \DEF{objective scalar equation}
if its weak form \EQU{*.weak}
for observer transformations $y=Y(y^*)$ obeys the transformation rule
\begin{equ}{eta}
  \eta^*=\eta\circ Y
\end{equ}
for scalar valued test functions $\eta\in\C_0^\infty(\UU;\RR)$.
That means that for $\eta^*$ as in \EQU{eta} and with $\UU=Y(\UU^*)$
\[
  \int_\UU\big( 
    \tx{\grad}\eta\dd q + \eta\rate
  \big)\d\leb4
  =
  \int_{\UU^*}\big( 
    \tx{\grad}\eta^*\dd q^* + \eta^*\rate^*
  \big)\d\leb4
\,.\]
This has been proved by a general theorem
(see \REF{continuum.thm.})
and it is true
if the quantities $q$ and $\rate$ satisfy the transformation rules
\begin{equ}{o}\begin{arr}l
  q\circ Y = \D{Y}q^*
\,,\quad
  \rate\circ Y = \rate^*
\,,\end{arr}\end{equ}
which written in components are
\[\begin{arr}{ll}
  q_k\circ Y = \sum_{l=0}^{n}Y_{k\p l}q^*_l
\,,\quad
  \rate\circ Y = \rate^*
\,.\end{arr}\]
Quantities with this property, that is, $q$ is an contravariant vector and $\rate$ an objective scalar,
are well-known.
Let us now present some special classes of scalar equations. 

\TITEL{Distributional form}
A more general version is
the distributional version
\begin{equ}{distr}
  \tx{\div}\,q=\rate \tm{in} \dis'(\UU;\RR)
\,,\end{equ}
where $q\in\dis'(\UU;\RR^n)$ und $\rate\in\dis'(\UU;\RR)$
are distributions. This definition means
for test functions $\eta\in\dis(\UU;\RR^n)=\C_0^\infty(\UU;\RR)$
\begin{equ}{distr-weak}
  \dprod{\tx\grad\eta}{q}{\dis(\UU)}+\dprod{\eta}{\rate}{\dis(\UU)}=0
\,,\end{equ}
and for this the transformation rule \EQU{eta} applies.
We add this distributional version because we will use it for
mass points in section \REF{point..}.
Similarly the distributional setting for other 
differential equations is defined.

\TITEL{Mass equation}
The standard example of a scalar equation is the mass equation. In this case the flux $q$ has the form
\begin{equ}{rhov}
  q=\rho\tx v+\JJ
\,.\end{equ}
We declare the \DEF{mass equation} as an objective scalar equation
\begin{equ}{mass}
  \tx\div(\rho\tx v+\JJ)=\rate
\,,\end{equ}
where the type of $\rho$ and $\tx v$ are defined in \REF{def-mass.}
and \REF{def-velo.} below so that $\rho\tx v$ is a contravariant vector.
Also $\JJ$ is assumed to be a contravariant vector
and $\rate$ an objective scalar.
Thus the condition \EQU{o} is satisfied.

\begin{stmt}{def-mass}{Mass density (Definition)}
A quantity $\rho\geq0$ is called \DEF{mass density},
if it is an objective scalar, that is, if for different observers
$
  \rho\circ Y=\rho^*
$.
A mass density is a density over space and time (this is essential).
\end{stmt}
This means that if $\rho$ is continuous in space and time then \begin{Equ}{density}
  \rho(y)=\lim_{U(y)\to\{y\}}\frac1{\leb4(U(y))}\int_{U(y)}\rho\d\leb4
\,,\end{Equ}
where $U(y)$ is a neighbourhood of the spacetime point $y$.
And for another observer with an observer transformation $y=Y(y^*)$
the \DEF{mass in $U(y)$} is
\[
  \int_{U(y)}\rho\d\leb4=\int_{U^*(y^*)}\rho^*\d\leb4
  \quad\tm{if $U(y)=Y(U^*(y^*))$}
\,,\]
where $\leb4(U(y))=\leb4(U^*(y^*))$
because we assume that $\fcn{det}\tx\D{Y}=1$.
(If $\rho$ is only integrable \EQU{density} holds only almost everywhere
in space and time.)

\begin{stmt}{def-velo}{Velocity (Definition)}
A quantity $\tx v$ is called a \DEF{4-velocity}, if it is an objective 4-vector, that is,
if for different observers
$
  \tx v\circ Y=\tx\D{Y}\tx v^*
$,
the rule for a contravariant vector,
and if with the (dual) time vector $\ewelt$ in section \REF{time..}
(see \EQU{time.*.welt})
\begin{equ}{o}
  \ewelt\dd\tx v=1
\,.\end{equ}
This definition is objective, because $\ewelt$ is a covariant vector.
\end{stmt}
\begin{prf}{\DE{}\EN{of objectivity}} It is
\[\begin{arr}l
  (\ewelt\dd\tx v)\circ Y=\ewelt\circ Y\dd(\tx\D{Y}\tx v^*)
  =(\up{\tx\D{Y}}T\ewelt\circ Y)\dd\tx v^*=\ewelt^*\dd\tx v^*
\,,\end{arr}\]
since $\ewelt$ is a covariant and $\tx v$ a contravariant vector.
\end{prf}

Here a remark on $\JJ$ and $\rho$.

\begin{stmt}{J}{Lemma}
Let $q=\rho\tx v+\JJ$.
An often used condition is $\ewelt\dd\JJ=0$.
This condition is objective and implies that
$\rho=\ewelt\dd q$.
That is, if the condition on $\JJ$ is satisfied
the mass density $\rho$ is the $\ewelt$-component
of the 4-flux in the mass equation.
\end{stmt}
\begin{prf}{}
Since
$
  (\ewelt\dd\JJ)\circ Y=(\ewelt\circ Y)\dd\D{Y}\JJ^*
  =(\up{\D{Y}}T\ewelt\circ Y)\dd\JJ^*=\ewelt^*\dd\JJ^*
$,
the condition on $\ewelt\dd\JJ$ is objective.
And $\ewelt\dd q=\rho\ewelt\dd\tx v+\ewelt\JJ=\rho$.
\end{prf}

The word ``4-velocity'' in connection to literature is described in \REF{u.}, 
where we use the definition
\begin{Equ}{norm}
  ||\tx w||:=\sqrt{\sum_{i\geq1}|e'_i\dd\tx w|^2}
  \quad\tm{for}\tx w\in\RR^4
\,.\end{Equ}

\begin{stmt}{u}{Remark on velocity}
In literature a ``4-velocity'' is a contravariant vector $\tx u$ satisfying
\begin{equ}{G-1}
  {\tx u}\dd\up\GG{-1}{\tx u}=-\cc^2
\,.\end{equ}
\begin{rem}{Reference}
Usually one has a different coordinate system, but it says that
$\big(\frac1{\cc}{\tx u}\big)\dd\tilde\GG\big(\frac1{\cc}{\tx u}\big)=1$
with $\tilde\GG:=-{\up\GG{-1}}$
(see for example \CIT{LandauLifschitz2}{I\,{\S}7}).
\end{rem}
\begin{enum}
\num{o}
Equation \EQU{G-1} is objective, that is, is the same for all observers.
\num{u} If $\tx u$ satisfies \EQU{G-1} then
$\ewelt\dd\tx u\neq0$ and (if $\ewelt\dd\tx u>0$)
  \[\tx v:=\frac{\tx u}{\ewelt\dd\tx u}\]
defines a 4-velocity as in \REF{def-velo.}.
It is $||\tx v||<\cc$.
\num{v} Let $\tx v$ be a 4-velocity as in \REF{def-velo.}
with $||\tx v||<\cc$ then
\[
  \tx u:=\gamma_{\tx v}\tx v\,,\quad
  \gamma_{\tx v}:=\frac1{\sqrt{1-\frac{||\tx v||^2}{\cc^2}}}
\,,\]
satisfies equation \EQU{G-1}.
\num{c} The classical formulas follow for $\GG=\GG_\cc$ and $\ewelt=\ee_0$.
\end{enum}
\end{stmt}
\begin{prf}{\REF{c}}
In the special case $\GG=\GG_\cc$ and $\ewelt=\ee_0$
(and hence $e'_k=e_k=\ee_k$ for $k\geq0$ in the standard case)
it follows that $\tx v=(1,v)$ and $||\tx v||=|v|$.
\end{prf}
\begin{prf}{\REF{o}} That ${\tx u}$ is contravariant
means ${\tx u}\circ Y=\D{Y}{\tx u^{*}}$, hence
\[\begin{arr}l
  ({\tx u}\dd\up\GG{-1}{\tx u})\circ Y
  =(\tx\D{Y}{\tx u^{*}})\dd\up\GG{-1}\tx\D{Y}{\tx u^{*}}
\\
  ={\tx u^{*}}\dd\big(
   \up{\tx\D{Y}}T\up\GG{-1}\tx\D{Y}
  \big){\tx u^{*}}
  ={\tx u^{*}}\dd\up{(\GG^*)}{-1}{\tx u^{*}}
\,,\end{arr}\]
since $\up\GG{-1}$ is a covariant tensor.
\end{prf}
\begin{prf}{\REF{u}} It is $\ewelt=e_0'$ and
\[
  \up\GG{-1}=-\cc^2e'_0\up{e'_0}T+\sum_{i\geq1}e'_i\up{e'_i}T
\,.\]
Hence
\[
  -\cc^2={\tx u}\dd\up\GG{-1}{\tx u}
  =  -\cc^2|{e'_0\dd\tx u}|^2+\sum_{i\geq1}|e'_i\dd\tx u|^2
\]
or
\[
  |{e'_0\dd\tx u}|^2=1+\frac1{\cc^2}\sum_{i\geq1}|e'_i\dd\tx u|^2
  =1+\frac{||\tx u||^2}{\cc^2}
\]
hence $|{e'_0\dd\tx u}|\geq1>0$.
If ${e'_0\dd\tx u}>0$ this means
\[
  {e'_0\dd\tx u}  =\sqrt{1+\frac{||\tx u||^2}{\cc^2}}
\,.\]
Then, with $\tx u_k:=e'_k\dd\tx u$ for $k\geq0$ we can write
$
  \tx u=\sum_{k\geq0}\tx u_ke_k
$
and we obtain
\[
  \tx v:=\frac1{e'_0\dd\tx u}\tx u=\frac{\tx u}{\tx u_0}
  =e_0+\sum_{i\geq1}\frac{\tx u_i}{\tx u_0}e_i
\,.\]
From this it follows immediately that $e'_0\dd\tx v=1$ and
\[
  ||\tx v||^2=\sum_{i\geq1}\Big|\frac{\tx u_i}{\tx u_0}\Big|^2
  =\frac{||\tx u||^2}{1+\frac{||\tx u||^2}{\cc^2}}<\cc^2
\,.\]
\end{prf}
\begin{prf}{\REF{v}} Let $\tx u=\mu\tx v$. It should be
\[
  -\cc^2={\tx u}\dd\up\GG{-1}{\tx u}
  =\mu^2(-\cc^2|\tx v_0|^2+\sum_{i\geq1}|\tx v_i|^2)
\]
where $\tx v_k:=e'_k\dd\tx v$. Since $\tx v_0=1$ and since
$||\tx v||<\cc$ this means
\[
  \mu^2=\left(1-\frac{||\tx v||^2}{\cc^2}\right)^{-1}
.\]
This proves the assertion.
\end{prf}

\TITEL{Diffusion}
A mass equation \EQU{mass} with a nonzero 4-flux $\JJ$ is called
\DEF{diffusion}.
Usually $\JJ$ depends on the gradient of an objective scalar $\mu$,
for example, $\mu$ is the mass density or the chemical potential.
\begin{stmt}{diffusion}{Theorem}
Let $\mu$ be an objective scalar. Then the following is true.
\begin{enum}
\num{Jfull} $\tx\JJ:=-\GG\tx\grad\mu$ is a contravariant vector.
\num{Jraum} $\JJ=-\GG^\raum\tx\grad\mu$ is a contravariant vector.
\num{eJraum} $\JJ$ satisfies the objective property
$\ewelt\dd\JJ=0$.
\num{eJfull} $\tx\JJ$ satisfies the objective property
$\ewelt\dd\tx\JJ=-(\GG\ewelt)\dd\tx\grad\mu=\frac1{\cc^2}\de_{e_0}\mu$.
\end{enum}
\end{stmt}

The diffusion equation becomes in the case \REF{diffusion.Jfull}
\begin{equ}{diffusion}
  \tx\div(\rho\tx v-\GG\tx\grad\mu)=\rate
\,.\end{equ}
Before we prove this let us remark the following.
In general we have the identities (see \EQU{time.thm.o})
\[
  \GG=-\frac1{\cc^2} e_0\up{e_0}T+\sum_{i=1}^ne_i\up{e_i}T
\,,\quad
  \ewelt=e_0'
\,.\]
\DE{}
\EN{By default $\GG$ and $\ewelt$ belong to those quantities
appearing in the description of physical processes.}
Now consider the splitting
\[
  \GG=\GG^\zeit+\GG^\raum
\,,\quad \GG^\zeit:=-\frac1{\cc^2}e_0\up{e_0}T=(\GG e_0')\up{e_0}T
\,,\]
that is, we get the identity \[
  \GG^\raum=\GG-\GG^\zeit=\GG\big(\tx{\Id}-e_0'\up{e_0}T\big)
\]
and this contains also $e_0$.
On the other hand $\ewelt\dd\tx\JJ=-(\GG\ewelt)\dd\tx\grad\mu$
only contains $\GG\ewelt=-\frac1{\cc^2}e_0$, which
in the classical limit $\cc\to\infty$ goes to zero.
But concerning $\JJ$ we also have $\GG^\raum\to\GG_\infty$ in the classical limit.
So one has to clarify what in nature the diffusion does.

\begin{prf}{\REF{Jfull}} Since $\mu$ is an objective scalar, that is
$\mu\circ Y=\mu^*$, we compute
\[
  \de_{\bar i}\mu^*=\de_{\bar i}(\mu\circ Y)
  =\sum_i(\de_i\mu)\circ Y\,Y_{i\p\bar i}
\]
for $\bar i=0\til 4$, hence
\[
  \tx\grad_{y^*}\mu^*=\up{\tx\D{Y}}T\tx\grad_{y}\mu
\,,\]
that is $\tx\grad\mu$ is a covariant vector.
Then we obtain
\[\begin{arr}l
  -\tx\JJ\circ Y=(\GG\circ Y)(\tx\grad\mu\circ Y)
  =\tx\D{Y}\GG^*\up{\tx\D{Y}}T\tx\grad\mu\circ Y
\\
  =\tx\D{Y}\GG^*\tx\grad\mu^*=-\tx\D{Y}\tx\JJ^*
\,,\end{arr}\]
that is, $\tx\JJ$ is an objective 4-vector
or a contravariant vector.
\end{prf}
\begin{prf}{\REF{eJraum}}We compute 
\[\begin{arr}l
  \tx\JJ-(e_0'\dd\tx\JJ)\,e_0
  =-\GG\tx\grad\mu+(e_0'\dd\GG\tx\grad\mu)\,e_0
  =-\GG\tx\grad\mu+\big((\GG e_0')\dd\tx\grad\mu\big)\,e_0
\\
  =-\GG\tx\grad\mu-\frac1{\cc^2}(e_0\dd\tx\grad\mu)\,e_0
  =-\Big(\GG+\frac1{\cc^2}e_0\up{e_0}T\Big)\tx\grad\mu
  =-\GG^\raum\tx\grad\mu=\JJ
\,,\end{arr}\]
hence $\ewelt\dd\JJ=0$.
\end{prf}

\begin{EXERCISES}
For later use the following \begin{stmt}{relativistic}{Example}Let $G$ be symmetric. Define the 4-flux $\JJ$ by
\[
  \JJ_i=\sum_{k}w_kG_{ki} \quad\tm{for}i\geq0
\,,\]
where $w$ is a covariant vector and $G$ a contravariant tensor.
\begin{enum}
\num{0} Then $\JJ$ is a contravariant vector.
\num{1} If $\ewelt\dd(Gw)=0$ then $\ewelt\dd\JJ=0$.
\num{2} If $G=\GG^\raum$ then $\ewelt\dd\JJ=0$.
\end{enum}
\end{stmt}
\begin{prf}{\REF{0}} Since $G$ is a contravariant tensor, that is
$G\circ Y=\tx DYG^*\up{\tx DY}T$, we conclude
\[\begin{arr}l
  \JJ_i\circ Y=\sum_{k}(w_kG_{ki})\circ Y
  =\sum_{k\bar k\bar i}w_k\circ Y\,Y_{k\p\bar k}Y_{i\p\bar i}G^*_{\bar k\bar i}
\\
  =\sum_{\bar i}Y_{i\p\bar i}\sum_{k\bar k}Y_{k\p\bar k}w_k\circ Y\,G^*_{\bar k\bar i}
  =\sum_{\bar i}Y_{i\p\bar i}\sum_{\bar k}w^*_{\bar k}G^*_{\bar k\bar i}
  =\sum_{\bar i}Y_{i\p\bar i}\JJ^*_{\bar i}
\,,\end{arr}\]
da $w^*_{\bar k}=\sum_{k}Y_{k\p\bar k}w_k\circ Y$.
\end{prf}
\begin{prf}{\REF{1}}
$\ewelt\dd\JJ=\ewelt\dd(\up{G}Tw)=\ewelt\dd(Gw)=0$.
\end{prf}
\begin{prf}{\REF{2}}
$\ewelt\dd\JJ=\up{\JJ}T\ewelt=\up{(\up{G}Tw)}T\ewelt=\up{w}TG\ewelt
=\up{w}T\GG^\raum\ewelt=0$.
\end{prf}

We shall later in \REF{fluid.J.} have a reason for such $\JJ$.
\end{EXERCISES}

{}}
\newpage
{}
\sect{impuls}{\DE{}\EN{Momentum equation}}
We introduce here the form of a relativistic momentum equation,
that is in a domain $\UU\subset\RR^4$ we look at the differential equation
\begin{equ}{strong}
  \tx{\div}\,T=\tx r
\end{equ}
for a 4-tensor $T=\seq{T_{ij}}{i,j=0\til3}$ and
a 4-field $\tx r=\seq{\tx r_{i}}{i=0\til3}$ which contain the 4-forces.
Written in coordinates this is
\begin{Equ}{sys}
  \sum_{j=0}^3\de_{y_j}T_{ij}=\tx r_i
  \quad\tm{for} i=0\til3
\,.\end{Equ}
It is essential to introduce
the weak formulation of this law and it reads for test functions $\zeta\in\C_0^\infty(\UU;\RR^4)$
\begin{Equ}{weak}
  \int_\UU\big( 
    \sum_{i,j=0}^3\de_{y_j}\zeta_i\cdot T_{ij} + \sum_{i=0}^3\zeta_i\,\tx r_i
  \big)\d\leb4=0
\,.\end{Equ}
The system \EQU{strong} is called \DEF{4-momentum equation}
if under observer transformations $Y$ the test functions
in \EQU{weak} transform via \begin{equ}{zeta}
  \zeta^*=\up{(\tx\D{Y})}T\zeta\circ Y
\,.\end{equ}
Here we as always assume that $\fcn{det}\tx\D{Y}=1$. This means with $\UU=Y(\UU^*)$ \[\begin{arr}l
  \int_\UU\big( 
    \sum_{i,j=0}^3\de_{y_j}\zeta_i\cdot T_{ij} + \sum_{i=0}^3\zeta_i\tx r_i
  \big)\d\leb4
  =
  \int_{\UU^*}\big( 
    \sum_{\bar i,\bar j=0}^3\de_{y^*_{\bar j}}\zeta^*_{\bar i}\cdot T^*_{\bar i\bar j}
    + \sum_{\bar i=0}^3\zeta^*_{\bar i}\tx r^*_{\bar i}
  \big)\d\leb4
\,.\end{arr}\]
Here the test function $\zeta$ is a covariant vector by \EQU{zeta}.
This definition is satisfied (see \REF{continuum.thm.})
if the quantities $T$ and $\tx r$ satisfy the transformation rules
\begin{Equ}{rule}\begin{arr}l
  T_{ij}\circ Y = \sum_{\bar i,\bar j=0}^3Y_{i\p\bar i}Y_{j\p\bar j}T^*_{\bar i\bar j}
\,,\\
  \tx r_i\circ Y = 
  \sum_{\bar i,\bar j=0}^3Y_{i\p\bar i\bar j}T^*_{\bar i\bar j}
  +\sum_{\bar i=0}^3Y_{i\p\bar i}\,\tx r^*_{\bar i}
\end{arr}\end{Equ}
for $i,j=0\til3$.
We mention that it is useful to make the following definition concerning the
transformation rule \EQU{rule}. Since $\tx r\circ Y$
depends linearly on $T^*$ we can write $\tx r$ as in \EQU{coriolis.+}.
After doing so the 4-momentum equation \EQU{sys} reads
\begin{Equ}{sys-C}
  \sum_{j=0}^3\de_{y_j}T_{ij}-\sum_{p,q=0}^3\coriolis^{pq}_iT_{pq}=\tx\force_i
  \quad\tm{for} i=0\til3
\,,\end{Equ}
where now $\tx\force$ is the real \DEF{4-force}
besides the fictitious forces containing forces
which we call \DEF{Coriolis forces}
(containing acceleration). The Coriolis coefficients $\coriolis^{pq}_i$ are defined in \REF{coriolis.}
and this definition takes care
of the terms with second derivatives $Y_{\p\bar i\bar j}$ in the
transformation rule \EQU{rule}.

\begin{stmt}{coriolis}{Coriolis Coefficients}
With $\coriolis_i:=\seq{\coriolis^{pq}_i}{p,q\geq0}$ let us write
\begin{Equ}{+}
  \tx r_i=\tx\force_i+\sum_{p,q\geq0}\coriolis^{pq}_iT_{pq}
\quad
  \text{for $i\geq0$,}
\end{Equ}
where $\coriolis^{pq}_i=\coriolis^{qp}_i$.
Then the rule \EQU{*.rule} for $\tx r$ is equivalent to
$\tx\force\circ Y=\tx\D{Y}\tx\force^*$ and
\begin{Equ}{o}
  \sum_{p,q\geq0}Y_{p\p\bar p}Y_{q\p\bar q}\coriolis^{pq}_i\circ Y
  =\sum_{\bar i\geq0}Y_{i\p\bar i}\coriolis^{*\bar p\bar q}_{\bar i} + Y_{i\p\bar p\bar q}
\end{Equ}
for all $i$ and $(\bar p,\bar q)$.
\end{stmt}
The matrix version of the above transformation rule \EQU{coriolis.o} is
\begin{Equ}{o-matrix}
  \up{\tx\D{Y}}T\coriolis_i\circ Y\,\tx\D{Y}=
  \sum_{\bar i\geq0}Y_{i\p\bar i}\coriolis^*_{\bar i}+\tx\D^2{Y_i}
.\end{Equ}
\begin{prf}{}
Using the above transformation rule for $\tx\force$
the second equation of \EQU{*.rule} becomes
\[
  \sum_{p,q}\coriolis^{pq}_i\circ Y\,T_{pq}\circ Y
  = \sum_{\bar p,\bar q}Y_{i\p\bar p\bar q}T^*_{\bar p\bar q}
  +\sum_{\bar i}Y_{i\p\bar i}
  \sum_{\bar p,\bar q}\coriolis^{*\bar p\bar q}_{\bar i}T^*_{\bar p\bar q}
\,.\]
Inserting on the left side for $T\circ Y$
the first equation of \EQU{*.rule} gives
\[
  \sum_{\bar p,\bar q}\sum_{p,q}\coriolis^{pq}_i\circ Y\,
  Y_{p\p\bar p}Y_{q\p\bar q}T^*_{\bar p\bar q}
  = \sum_{\bar p,\bar q}Y_{i\p\bar p\bar q}T^*_{\bar p\bar q}
  +\sum_{\bar i}Y_{i\p\bar i}
  \sum_{\bar p,\bar q}\coriolis^{*\bar p\bar q}_{\bar i}T^*_{\bar p\bar q}
\,.\]
Now compare the coefficients of $T^*_{\bar p\bar q}$.
\end{prf}

\begin{stmt}{klass}{\DE{}\EN{Remark (Classical physics)}}
\DE{}
\EN{In the classical limit $\cc\to\infty$ the derivatives $Y_{k\p x_ix_j}$
are zero, only the second derivatives}
\[
\left.\begin{arrc}l
  Y_{k\p tt}
\\
  Y_{k\p tx_j}=Y_{k\p x_jt}
\end{arrc}\right\}
\tm{for $j,k\geq1$ can be nonzero.}
\]
Or, if we write $Y_k=X_k$ for $k\geq1$, 
where $X(t^*,x^*)=Q(t^*)x^*+\xdiff(t^*)$,
the nonzero terms are those depending on $\ddot X$ and $\dot Q$.
Therefore the Coriolis coefficient $\coriolis$ is of the form
\[
  \coriolis_0=0 \,,\quad
  \coriolis_k=\begin{mat}a_k\tab\up{d_k}T\\d_k\tab0\end{mat}
  \tm{for}k\geq1
\]
with
\begin{equ}{+}\begin{arr}c
  \up{Q}Td_k = \sum_{\bar k}Q_{k\bar k}d^*_{\bar k} + \seq{\dot Q_{kj}}{j}
\,,\\
  a_k+2\dot X\dd d_k = \sum_{\bar k}Q_{k\bar k}a^*_{\bar k} + \ddot X_k
\,.\end{arr}\end{equ}
This gives for the Coriolis terms the well known fictitious forces
\[
  \sum_{p,q\geq0}\coriolis^{pq}_kT_{pq}
  =a_kT_{00}+\sum_{p\geq0}d_{kp}(T_{k0}+T_{0k})
  =a_k\rho+\sum_{p\geq0}d_{kp}(2\rho v_k+\JJ_k)
\,.\]
\end{stmt}
\begin{rem}{Remark} In the classical limit
this has been computed with different notation
in \CIT{Alt-FluidSolid}{9 Force}. \end{rem}
\begin{prf}{}
If $Y$ is a Newton transformation, then the inhomogeneous term in
the transformation rule for $\tx r$ is for $k\geq1$
\[
  \sum_{\bar i,\bar j\geq0}Y_{k\p\bar i\bar j}T^*_{\bar i\bar j}
  =\ddot X_kT^*_{00}+\sum_{\bar j\geq1}\dot Q_{kj}(T^*_{0\bar j}+T^*_{\bar j0})
\,.\]
Here the first term is the $k^{\rm th}$-component of the acceleration
and the second term the $k^{\rm th}$-component of the Coriolis force.
The transformation rule for $\coriolis$ is
$\coriolis_0=0$ and for $k\geq1$
\[\begin{arr}l
  \begin{mat}1\tab\up{\dot X}T\\0\tab\up{Q}T\end{mat}
  \begin{mat}a_k\tab\up{d_k}T\\d_k\tab0\end{mat}
  \begin{mat}1\tab0\\\dot X\tab Q\end{mat}
\\
  =\sum_{\bar k}Q_{k\bar k}
  \begin{mat}a^*_{\bar k}\tab\up{d^*_{\bar k}}T\\d^*_{\bar k}\tab0\end{mat}
  +\begin{mat}\ddot X_k\tab\up{\seq{\dot Q_{kj}}{j}}T
  \\\seq{\dot Q_{kj}}{j}\tab0\end{mat}
,\end{arr}\]
which is equivalent to \EQU{+}.
These are the transformation rules for $\seq{a_k}{k}$ and $\seq{d_k}{k}$.
\end{prf}

The transformation rule \EQU{o-matrix} applies to
$\coriolis=\seq{\coriolis^{pq}_k}{kpq}$
and to the negative Christoffel symbols
$\christoffel=\seq{\christoffel_{pq}^k}{kpq}$,
see \REF{christoffel.}. Thus $B^{pq}_k:=\coriolis^{pq}_k+\christoffel_{pq}^k$ satisfy the rule
\[\up{\tx\D{Y}}TB_k\circ Y\,\tx\D{Y}=\sum_{\bar k\geq0}Y_{k\p\bar k}B^*_{\bar k}\,.\]
Therefore if the $B_k$ vanish for one ob\-ser\-ver, it is always true that
$\coriolis^{pq}_k=-\christoffel_{pq}^k$.

\begin{stmt}{christoffel}{Christoffel symbols} They are
\[
  \christoffel_{ij}^k:=
  \frac12\sum_lg^{kl}\big(
  g_{jl\p i}+g_{il\p j}-g_{ij\p l}
  \big)
\]
for $i,j,k\geq0$ and $-\christoffel$ satisfies the transformation rule
\EQU{coriolis.o} (formulated for the Coriolis coefficients). Here we use
\begin{equ}{def}
  g^{kl}=-\GG_{kl} \quad\tm{and}\quad g_{kl}=-{\up\GG{-1}}_{kl}
\,.\end{equ}
\end{stmt}
\begin{APPEND}
\begin{prf}{}
This is well known, see \CIT{LandauLifschitz2}{{\S}86}.
We have the following transformation rules 
\[\begin{arr}l
  g^{kl}\circ Y=\sum_{nm}Y_{k\p n}Y_{l\p m}\,g^{*nm}
  \tm{for}y=Y(y^*)
\,,\\
  g_{ij}=\sum_{pq}Y^*_{p\p i}Y^*_{q\p j}\,g^*_{pq}\circ Y^*
  \tm{for}y^*=Y^*(y)
\,.\end{arr}\]
Now inserting this, we obtain
\[\begin{arr}l
  \Gamma^k_{ij}=\frac12\sum_{l}g^{kl}\big( g_{jl\p i}+g_{il\p j}-g_{ij\p l} \big)
\\
  =\frac12\sum_{lnmpq}(g^{*nm}Y_{k\p n}Y_{l\p m})\circ Y^*\Big( 
  \big(g^*_{pq}\circ Y^* Y^*_{p\p j}Y^*_{q\p l}\big)_{\p i}
\\\hfill
  +\big(g^*_{pq}\circ Y^* Y^*_{p\p i}Y^*_{q\p l}\big)_{\p j}
  -\big(g^*_{pq}\circ Y^* Y^*_{p\p i}Y^*_{q\p j}\big)_{\p l}
  \Big)
\\
  =\frac12\sum_{lnmpq}(g^{*nm}Y_{k\p n}Y_{l\p m})\circ Y^*\Big( 
  \big(g^*_{pq}\circ Y^*\big)_{\p i} Y^*_{p\p j}Y^*_{q\p l}
\\\hfill
  +\big(g^*_{pq}\circ Y^*\big)_{\p j} Y^*_{p\p i}Y^*_{q\p l}
  -\big(g^*_{pq}\circ Y^*\big)_{\p l} Y^*_{p\p i}Y^*_{q\p j}
  \Big)
\\
  +\frac12\sum_{lnmpq}(g^{*nm}Y_{k\p n}Y_{l\p m})\circ Y^*\Big( 
  g^*_{pq}\circ Y^* \big(Y^*_{p\p j}Y^*_{q\p l}\big)_{\p i}
\\\hfill
  +g^*_{pq}\circ Y^* \big(Y^*_{p\p i}Y^*_{q\p l}\big)_{\p j}
  -g^*_{pq}\circ Y^* \big(Y^*_{p\p i}Y^*_{q\p j}\big)_{\p l}
  \Big)
\,.\end{arr}\]
The first sum is
\newpage
\[\begin{arr}l
  \frac12\sum_{lnmpq}(g^{*nm}Y_{k\p n}Y_{l\p m})\circ Y^*\Big( 
  \big(g^*_{pq}\circ Y^*\big)_{\p i} Y^*_{p\p j}Y^*_{q\p l}
\\\hfill
  +\big(g^*_{pq}\circ Y^*\big)_{\p j} Y^*_{p\p i}Y^*_{q\p l}
  -\big(g^*_{pq}\circ Y^*\big)_{\p l} Y^*_{p\p i}Y^*_{q\p j}
  \Big)
\\
  =\frac12\sum_{nmpq}(g^{*nm}Y_{k\p n})\circ Y^*\Big(
  \big(g^*_{pq}\circ Y^*\big)_{\p i}Y^*_{p\p j}
  \Ubrack{\sum_lY^*_{q\p l}Y_{l\p m}\circ Y^*}{=\kronecker{qm}}
\\\hfill
  +\big(g^*_{pq}\circ Y^*\big)_{\p j}Y^*_{p\p i}
  \Ubrack{\sum_lY^*_{q\p l}Y_{l\p m}\circ Y^*}{=\kronecker{qm}}
  -\sum_l\big(g^*_{pq}\circ Y^*\big)_{\p l}Y_{l\p m}\circ Y^*\,Y^*_{p\p i}Y^*_{q\p j}
  \Big)
\\
  =\frac12\sum_{nmp}(g^{*nm}Y_{k\p n})\circ Y^*\Big( 
  \big(g^*_{pm}\circ Y^*\big)_{\p i}Y^*_{p\p j}
  +\big(g^*_{pm}\circ Y^*\big)_{\p j}Y^*_{p\p i}
\\\hfill
  -\sum_{qr}g^*_{pq\p r}\circ Y^*
  \Ubrack{\sum_lY^*_{r\p l}Y_{l\p m}\circ Y^*}{=\kronecker{rm}}
  Y^*_{p\p i}Y^*_{q\p j}
  \Big)
\\
  =\frac12\sum_{nm}(g^{*nm}Y_{k\p n})\circ Y^*\Big( 
  \Ubrack{\sum_{pq}g^*_{pm\p q}Y^*_{q\p i}Y^*_{p\p j}}
  {=\sum_{pq}g^*_{qm\p p}Y^*_{p\p i}Y^*_{q\p j}}
  +\sum_{pq}g^*_{pm\p q}Y^*_{q\p j}Y^*_{p\p i}
\\\hfill
  -\sum_{pq}g^*_{pq\p m}\circ Y^*\,Y^*_{p\p i}Y^*_{q\p j}
  \Big)
\\
  =\frac12\sum_{nmpq}(g^{*nm}Y_{k\p n})\circ Y^*\big(
  g^*_{qm\p p}+ g^*_{pm\p q}-g^*_{pq\p m}
  \big)\circ Y^*\,Y^*_{p\p i}Y^*_{q\p j}
\,.\end{arr}\]
The second sum is
\[\begin{arr}l
  \frac12\sum_{lnmpq}(g^{*nm}Y_{k\p n}Y_{l\p m})\circ Y^*\Big( 
  g^*_{pq}\circ Y^* \big(Y^*_{p\p j}Y^*_{q\p l}\big)_{\p i}
\\\hfill
  +g^*_{pq}\circ Y^* \big(Y^*_{p\p i}Y^*_{q\p l}\big)_{\p j}
  -g^*_{pq}\circ Y^* \big(Y^*_{p\p i}Y^*_{q\p j}\big)_{\p l}
  \Big)
\\
  =\frac12\sum_{lnmpq}(g^{*nm}g^*_{pq}Y_{k\p n}Y_{l\p m})\circ Y^*\Big( 
  Y^*_{p\p ij}Y^*_{q\p l}+Y^*_{p\p j}Y^*_{q\p il}
\\\hfill
  +Y^*_{p\p ji}Y^*_{q\p l}+Y^*_{p\p i}Y^*_{q\p jl}
  -Y^*_{p\p il}Y^*_{q\p j}-Y^*_{p\p i}Y^*_{q\p jl}
  \Big)
\\
  =\sum_{lnmpq}(g^{*nm}g^*_{pq}Y_{k\p n}Y_{l\p m})\circ Y^*\,Y^*_{p\p ij}Y^*_{q\p l}
\\
  =\sum_{nmpq}(g^{*nm}g^*_{pq}Y_{k\p n})\circ Y^*\,Y^*_{p\p ij}
  \sum_lY^*_{q\p l}Y_{l\p m}\circ Y^*
\\
  =\sum_{nmp}(g^{*nm}g^*_{pm}Y_{k\p n})\circ Y^*\,Y^*_{p\p ij}
  =\sum_{n}Y_{k\p n}\circ Y^*\,Y^*_{n\p ij}
\,.\end{arr}\]
Together we have shown that
\[\begin{arr}l
  \Gamma^k_{ij}=\frac12\sum_{l}g^{kl}\big( g_{jl\p i}+g_{il\p j}-g_{ij\p l} \big)
\\
  =\frac12\sum_{nmpq}(g^{*nm}Y_{k\p n})\circ Y^*\big(
  g^*_{qm\p p}+ g^*_{pm\p q}-g^*_{pq\p m}
  \big)\circ Y^*\,Y^*_{p\p i}Y^*_{q\p j}
\\\hfill
  +\sum_{n}Y_{k\p n}\circ Y^*\,Y^*_{n\p ij}
\\
  =\sum_{n}Y_{k\p n}\circ Y^*\Big(
  \frac12\sum_{mpq}\big(g^{*nm}(g^*_{qm\p p}+ g^*_{pm\p q}-g^*_{pq\p m})
  \big)\circ Y^*\,Y^*_{p\p i}Y^*_{q\p j}
  +Y^*_{n\p ij}
  \Big)
\\
  =\sum_{n}Y_{k\p n}\circ Y^*\Big(
  \sum_{pq}\Gamma^{*n}_{pq}\circ Y^*\,Y^*_{p\p i}Y^*_{q\p j}
  +Y^*_{n\p ij}
  \Big)
\,,\end{arr}\]
that is
\begin{Equ}{star}\begin{arr}l
  \Gamma^k_{ij}
  =\sum_{n}Y_{k\p n}\circ Y^*\Big(
  \sum_{pq}\Gamma^{*n}_{pq}\circ Y^*\,Y^*_{p\p i}Y^*_{q\p j}
  +Y^*_{n\p ij}
  \Big)
\,.\end{arr}\end{Equ}
\DE{}
\EN{The equation \EQU{star} is equivalent to}
\[
  \sum_{k}Y^*_{n\p k}\Gamma^k_{ij}
  =\sum_{pq}\Gamma^{*n}_{pq}\circ Y^*\,Y^*_{p\p i}Y^*_{q\p j}
  +Y^*_{n\p ij}
\]
\DE{}
\EN{and this is equivalent to}
\[
  \sum_{ij}(Y_{i\p p}Y_{j\p q})\circ Y^*\,\Gamma^k_{ij}
  =\sum_{n}Y_{k\p n}\circ Y^*\,\Gamma^{*n}_{pq}\circ Y^*
  +\sum_{nij}(Y_{k\p n}Y_{i\p p}Y_{j\p q})\circ Y^*\,Y^*_{n\p ij}
\,.\]
\DE{}
\EN{Since}
\[
  -Y_{k\p pq}=
  \sum_{nij}(Y_{k\p n}Y_{i\p p}Y_{j\p q})\,Y^*_{n\p ij}\circ Y
\]
\DE{}
\EN{it is equivalent to}
\[
  \sum_{ij}Y_{i\p p}Y_{j\p q}\,\Gamma^k_{ij}\circ Y
  =\sum_{n}Y_{k\p n}\,\Gamma^{*n}_{pq}
  -Y_{k\p pq}
\]
\DE{}
\EN{which is \EQU{coriolis.o} for $-\Gamma$.}
\end{prf}
\end{APPEND}

Therefore it is of interest to take the Christoffel symbols
for a relativistic sequence and show that in the classical limit
they converge to
the formulas for the Coriolis coefficients in \REF{klass.}.
The following lemma holds.

\begin{stmt}{id}{\DE{}\EN{Identity for $\Gamma$}}
The Christoffel symbols satisfy (we use \EQU{christoffel.def})
\[
  g_{ik\p l}=\sum_m\big(
  g_{mk}\Gamma^m_{il}+g_{im}\Gamma^m_{kl}
  \big)
\,.\]
\end{stmt}
This is well known, see \CIT{LandauLifschitz2}{{\S}86}.

\medskip
\TITEL{Contains mass equation}
It is an essential step to prove that the 4-momentum system \EQU{sys} contains
as part the mass equation.
To explain this
let us assume for the moment that the 4-tensor $T$ has the form
\begin{equ}{vv}
  T_{ij}=\rho\tx v_i\tx v_j+\tx\Pi_{ij}
  \quad\tm{for}i,j\geq0
,\end{equ}
where $\tx v$ is a 4-velocity as in \REF{scalar.def-velo.}
and $\rho$ is an objective scalar, the mass density.
If we multiply this with the component $\ewelt_i$
of the vector $\ewelt$ and sum over $i$ we get
using that $\ewelt\dd\tx v=1$
\begin{equ}{qq}
  q_j:=\sum_i\ewelt_iT_{ij}=\rho\tx v_j+\sum_i\ewelt_i\tx\Pi_{ij}
  \quad\tm{for}j\geq0
,\end{equ}
where $q$ has the form \EQU{scalar.*.rhov}, that is,
$q$ is identical with a mass flux.
We convert this to an argument with test functions
by looking at the weak version \EQU{weak}.
By \EQU{zeta} the test function $\zeta$ must be a covariant vector.
With an objective scalar test function $\eta$ we can set $\zeta=\eta\ewelt$.
Since $\ewelt$ is a covariant vector
this is a possible test function. Thus we obtain a weak equation written in $\eta$
which is the ``$\ewelt$-component'' of the relativistic momentum equation \EQU{sys} and
is identical to the mass conservation \EQU{scalar.*.weak}.
Doing so for arbitrary $T$ we get the
\begin{stmt}{mass}{Mass conservation (Reduction)}
A part of the 4-momentum equation \EQU{*.strong} reads
\begin{equ}{o}
  \tx{\div}\,q=\rate
\end{equ}
with
\[\begin{arr}l
  q_j:=\sum_{i}\ewelt_iT_{ij}
\,,\quad
  \rate:=\sum_{i}\ewelt_i\tx r_i + \sum_{i,j}\ewelt_{i\p j}T_{ij}
\,.\end{arr}\]
If $T$ is as in \EQU{*.vv} then
\[\begin{arr}l
  q_j=\rho\tx v_j+\JJ_j
\,,\quad
  \JJ_j:=\sum_{i}\ewelt_i\tx\Pi_{ij}
\,,\\
  \rate
  \begin{arr}l
    :=\sum_{i}\ewelt_i\big(\tx r_i-\rho(\tx v\dd\tx\grad)\tx v_i\big)
    + \sum_{i,j}\ewelt_{i\p j}\tx\Pi_{ij}
  \\
    =\ewelt\dd\tx\force + \sum_{p,q}(\ewelt\dd\coriolis^{pq}+\ewelt_{p\p q})T_{pq}
  \end{arr}
\,.\end{arr}\]
\end{stmt}
\begin{rem}{Remark}
We mention that the condition $\sum_{i,j}\ewelt_i\ewelt_j\tx\Pi_{i,j}=0$
implies $\ewelt\dd\JJ=0$.
\end{rem}
\begin{prf}{}We take
\begin{equ}{reduce}
  \zeta:=\eta\,\ewelt
\end{equ}
as test function.
If $\eta$ is an objective scalar then
\[
  \zeta^*=\eta^*\ewelt^*=\eta^*\up{\D{Y}}T\ewelt\circ Y
  =\up{\D{Y}}T\zeta\circ Y
\,,\]
hence $\zeta$ is an allowed test function and
it follows from \EQU{*.weak}
\[\begin{arr}l
  0=\int_\UU\big( 
  \sum_{i,j}\de_{y_j}\zeta_i T_{ij} + \sum_{i}\zeta_i\tx r_i
  \big)\d\leb4
\\
  =\int_\UU\big(
  \sum_{i,j}\de_{y_j}(\eta\ewelt_i) T_{ij} + \eta\sum_{i}\ewelt_i\tx r_i
  \big)\d\leb4
\\
  =\int_\UU\Big(
  \sum_{j}\de_{y_j}\eta\Ubrack{\sum_{i}\ewelt_iT_{ij}}{=q_j}
  +\eta\big(
  \Ubrack{\sum_{i,j}(\de_{y_j}\ewelt_i)T_{ij} + \sum_{i}\ewelt_i\tx r_i}{=\rate}\big)
  \Big)\d\leb4
\,.\end{arr}\]
If $T$ is as in \EQU{*.vv} then
\[
  q_j=\sum_i\ewelt_i(\rho\tx v_i\tx v_j+\tx\Pi_{ij})
  =\rho(\ewelt\dd\tx v)\tx v_j+\sum_i\ewelt_i\Pi_{ij}
\]
with $\ewelt\dd\tx v=1$, and one term of
$\rate=\sum_{i}\ewelt_i\tx r_i + \sum_{i,j}\ewelt_{i\p j}T_{ij}$ is
\[
  \sum_{i,j}\ewelt_{i\p j}\rho\tx v_i\tx v_j
  =\sum_{i,j}(\ewelt_{i}\tx v_i)_{\p j}\rho\tx v_j
  -\sum_{i,j}\ewelt_{i}\tx v_{i\p j}\rho\tx v_j
  =-\rho\sum_{i}\ewelt_{i}\sum_{j}\tx v_j\tx v_{i\p j}
\]
again since $\ewelt\dd\tx v$ is constant,
and $(\tx v\dd\tx\grad)\tx v=\tx\D{\tx v}\,\tx v$. Also $\rate$ equals
\[\begin{arr}l
  \rate=\sum_{i}\ewelt_i\Big(
  \tx\force_i+\sum_{p,q}\coriolis_i^{pq}T_{pq}
  \Big) + \sum_{i,j}\ewelt_{i\p j}T_{ij}
\\
  =\sum_{i}\ewelt_i\tx\force_i
  +\sum_{p,q}\big(\ewelt_i\coriolis_i^{pq}+\ewelt_{p\p q}\big)T_{pq}
\,,\end{arr}\]
whish gives the second representation of $\rate$.
\end{prf}

It is now clear that the mass equation is a part of the
4-momentum system. A question is how one formulates in general
the momentum system without the mass equation.

\begin{stmt}{reduct}{Reduced momentum equation}
This equation reads
\[
  \int_{\RR^4}\big( 
  \sum_{i,j=0}^3\de_{y_j}\tilde\zeta_i\cdot T_{ij}
  +\sum_{i=0}^3\tilde\zeta_i\force_i
  \big)\d\leb4=0
  \quad\tm{for} \tilde\zeta\dd e_0=0
.\]
Here $\force$ is the remaining force (without reaction term)
\[
  \force:=\tx r-\tx r\dd e_0'e_0
\,.\]
\begin{rem}{Remark}
So this differential equation together with \EQU{mass.o}
is equivalent to the 4-momentum system \EQU{*.sys}.  
\end{rem}
\end{stmt}
\begin{prf}{}
The equation $\zeta=\eta e_0'$ implies $\eta=\zeta\dd e_0$. Or
$\tilde\zeta:=\zeta-\eta e_0'$ with $\eta:=\zeta\dd e_0$ gives
$\tilde\zeta\dd e_0=0$. And the definition of $\force$ gives
$\tilde\zeta\dd\tx r=\tilde\zeta\dd\force$.
Therefore
\[
  \zeta\dd\tx r=(\tilde\zeta+\eta e_0')\dd\tx r
  =\eta e_0'\dd\tx r + \tilde\zeta\dd\force
\,,\]
where $e_0'\dd\tx r=\rate-\rho e_0'\dd\tx\D{\tx v}\,\tx v$.
The $\eta$-term is the contribution in the mass equation
and the $\tilde\zeta$-term the contribution in the reduced system.
\end{prf}

{}
\newpage
{{}
\sect{point}{\DE{}\EN{Moving particle}}
The mass-momentum equation of section \REF{impuls..} should be consistent with
the law of moving particles.
Therefore we consider the differential equation \EQU{impuls.*.strong}
\[
  \tx\div\,T = \tx r
\]
now for distributions based on an evolving point $\Gamma\subset\RR^4$.
We refer to \REF{mass.} and \REF{momentum.} for the distributional version
of differential equations on $\Gamma$.

\medskip
By an \DEF{evolving point} $\Gamma\subset\RR^4$
we mean that for $y\in\Gamma$ the scalar product $\ewelt(y)\dd\tau>0$ 
for some $\tau\in \TT_y(\Gamma)\setminus\{0\}$,
where $\TT(\Gamma)$ denotes the spacetime tangent space of $\Gamma$.
Multiplying $\tau$ by a positive constant we obtain $\ewelt(y)\dd\tau=1$.
This way, since $\Gamma$ is one-dimensional, the tangential vector
becomes unique. Thus    we obtain a 4-vector
\begin{equ}{velo}
  \tx\velo_\Gamma(y)\in\TT_y(\Gamma)
\quad\text{with}\quad
  \ewelt(y)\dd\tx\velo_\Gamma(y)=1
,\end{equ}
which we call the \DEF{4-velocity} of $\Gamma$,
a notion which we already introduced in \REF{scalar.def-velo.}.

\begin{stmt}{4velocity}{Remark}
It follows from \EQU{*.velo} that $\tx\velo_\Gamma=e_0+\velo_\Gamma$ with
$\velo_\Gamma\in\welt=\up{\{\ewelt\}}\perp$.
And $\tx\velo_\Gamma$ is a contravariant vector,
i.e.~it satisfies the transformation rule
\[\tx\velo_\Gamma\circ Y=\D{Y}\tx\velo_{\Gamma^*}\,,\]
where $Y$ is the observer transformation.
\end{stmt}
\begin{prf}{}
If $\Gamma=Y(\Gamma^*)$ it follows for $y=Y(y^*)$
that $\tx\D{Y}(y^*)$
maps $\TT_{y^*}(\Gamma^*)$ into $\TT_y(\Gamma)$.
\end{prf}

In order to get an impression of a differential equation on $\Gamma$
we prove the following lemma.

\begin{stmt}{general}{Lemma (Objectivity on $\Gamma$)}
For an evolving point $\Gamma$ the following holds:
\[
  \int_\Gamma \frac{g}{|\tx\velo_{\Gamma}|}\d\haus1
  =\int_{\Gamma^*} \frac{g^*}{|\tx\velo_{\Gamma^*}|}\d\haus1
\,,\]
if $g\maps\Gamma\to\RR$ is an objective scalar.
The velocity $\tx\velo_\Gamma$ is the 4-velocity as in \REF{4velocity.}.
Here $|\tx\velo_\Gamma|$ denotes the Euclidean norm of $\tx\velo_\Gamma$.
\end{stmt}

\begin{prf}{}
Let $y=Y(y^*)$ be an observer transformation and $\Gamma=Y(\Gamma^*)$.
Then the transformation formula from mathematics says that for every local function
$f:\Gamma\to\RR$
\[
  \int_\Gamma f\d\haus1
  =\int_{\Gamma^*}f\circ Y
  \,|\restrict{\fcn{det}\D{Y}}{\TT_{y^*}(\Gamma^*)}{}|\d\haus1
\,.\]
We have to bring this in our statement. If we choose tangent vectors
\[
  \tau(y)=\D{Y}(y^*)\,\tau^*(y^*)\in\TT_{y}(\Gamma)
\quad\tm{\DE{}\EN{for}}\quad
  \tau^*(y^*)\in\TT_{y^*}(\Gamma^*)
\,,\]
and if we choose $\tau^*(y^*)$ as unit vector $|\tau^*(y^*)|=1$ then
\[
  |\restrict{\fcn{det}\D{Y}}{\TT_{y^*}(\Gamma^*)}{}|=|\tau\circ Y|
\,.\]
Since $\tx\velo_\Gamma(y)=\D{Y}(y^*)\tx\velo_{\Gamma^*}(y^*)$
by \REF{4velocity.} we can choose
\[
  \tau^*(y^*)=\frac{\tx\velo_{\Gamma^*}(y^*)}{|\tx\velo_{\Gamma^*}(y^*)|}
\,,\quad
  \tau(y)=\frac{\tx\velo_\Gamma(y)}{|\tx\velo_{\Gamma^*}(y^*)|}
\]
and obtain
\[
  |\restrict{\fcn{det}\D{Y}}{\TT_{y^*}(\Gamma^*)}{}|
  =\frac{|\tx\velo_\Gamma\circ Y|}{|\tx\velo_{\Gamma^*}|}
\,.\]
Therefore the above formula reads
\[
  \int_\Gamma f\d\haus1
  =\int_{\Gamma^*}(f\,|\tx\velo_\Gamma|)\circ Y\frac{\d\haus1}{|\tx\velo_{\Gamma^*}|}
\,.\]
Setting $g(y):=f(y)|\tx\velo_{\Gamma}(y)|$ we obtain
\[
  \int_\Gamma\frac{g}{|\tx\velo_{\Gamma}|}\d\haus1
  =\int_{\Gamma^*}\frac{g\circ Y}{|\tx\velo_{\Gamma^*}|}\d\haus1
\,.\]
\DE{}
\EN{If $g$ is an objective scalar, that is $g\circ Y=g^*$,
the assertion follows.}
\end{prf}

The 4-velocity $\tx\velo_{\Gamma}$ is a contravariant vector,
but the transformation rule for $|\tx\velo_{\Gamma}|$ is not so easy,
it is more convenient to consider the measure
\begin{Equ}{meas}
  \meas\Gamma:=\frac1{|\tx\velo_{\Gamma}|}\haus1\cutted\Gamma
\end{Equ}
and therefore one can write the result of \REF{general.} as
\begin{Equ}{int}
  \int_{\RR^4} g\d\meas{\Gamma}=\int_{\RR^4}g^*\d\meas{\Gamma^*}
\,.\end{Equ}
We can view this also as a transformation rule for $\meas{\Gamma}$:
\begin{stmt}{meas}{Remark}
If $Y$ is an observer transformation and $\Gamma=Y(\Gamma^*)$ then
\[
  \meas{\Gamma}(B)=\meas{\Gamma^*}(B^*)
  \tm{for} B=Y(B^*)
.\]
\end{stmt}

We now
write down differential equations on the curve $\Gamma$
in order to describe the movement of a ``tiny mass''.
First let us use the distributional version of the mass equation
\EQU{scalar.*.distr} taking \REF{general.} into account.

\begin{stmt}{mass}{Mass equation} We define the distributional mass equation by
\[\begin{arr}c
  \tx\div q=\rate
\,,\\
  {q}=m\tx v\meas\Gamma  \,,\quad  \rate=r\meas\Gamma
\,,\end{arr}\]
where the mass $m$ of the particle is an objective scalar
and the velocity $\tx v=\tx\velo_\Gamma$
and also the rate $r$ is an objective scalar.
For scalar test functions $\eta$ this reads
\[\begin{arr}l
  0=\dprod{\eta}{-\tx\div{q}+\rate}{\dis(\UU)}
  =\int_{\RR^4}\big(\tx\grad\eta\dd(m\tx v)+\eta r\big)\d\meas{\Gamma}
\\
  =\int_\Gamma\Big(\tx\grad\eta\dd(m\tx\velo_\Gamma)+\eta r\Big)
  \frac{\d\haus1}{|\tx\velo_\Gamma|}
  =\int_\Gamma\eta\Big(
  -\tx\div^\Gamma\Big(\frac{m\tx\velo_\Gamma}{|\tx\velo_\Gamma|}\Big)
  +\frac{r}{|\tx\velo_\Gamma|}
  \Big)\d\haus1
\end{arr}\]
hence
\[
  \tx\div^\Gamma\Big(\frac{m\tx\velo_\Gamma}{|\tx\velo_\Gamma|}\Big)
  =\frac{r}{|\tx\velo_\Gamma|}
\,.\]
For the last integral the regularity of $\Gamma$ and $m$ are required.

\begin{rem}{\DE{}\EN{Hint}}
We mention section~\REF{approximation..} where the notion of
the scalar $m$ is motivated.
\end{rem}
\end{stmt}
\begin{prf}{}
Since $\eta$ is an objective scalar, that is $\eta\circ Y=\eta^*$,
we deduce that $\grad\eta$ satisfies
$\tx\grad\eta^*=\up{(\tx\D{Y})}T\tx\grad\eta\circ Y$.
Then in the first integral the $\meas\Gamma$-integrand
\[
  g:=\tx\grad\eta\dd(m\tx v)+\eta r
\]
satisfies
\[\begin{arr}l
  g^*=\tx\grad\eta^*\dd(m^*\tx v^*)+\eta^*r^*
  =(\up{(\tx\D{Y})}T\tx\grad\eta\circ Y)\dd(m^*\tx v^*)
  +(\eta\circ Y)r^*
\\
  =(\tx\grad\eta\circ Y)\dd(m^*\tx\D{Y}v^*)+(\eta\circ Y)r^*
  =(\tx\grad\eta\dd(m\tx v)+\eta r)\circ Y=g\circ Y
\,,\end{arr}\]
that is, $g$ is an objective scalar.
Hence by \EQU{*.int}
\[
  \dprod{\eta}{-\div{q}+\rate}{\dis(\UU)}
  =\int_\Gamma g\d\meas{\Gamma}
  =\int_{\Gamma^*}g^*\d\meas{\Gamma^*}
  =\dprod{\eta^*}{-\div{q^*}+\rate^*}{\dis(\UU^*)}
\,.\]
This shows that the distributional equation is the mass equation
(see also section~\REF{approximation..}).
To derive the strong version note that
$m\tx v=m\tx\velo_\Gamma\in\TT(\Gamma)$.
\end{prf}

Similarly, we treat the momentum equation for a particle.

\begin{stmt}{momentum}{Momentum equation} The distributional 4-momentum system is defined by
\[\begin{arr}c
  \tx\div T=\tx r \quad\text{with}
\\
  T_{ij}=m\tx v_i\tx v_j\meas\Gamma  \,,\quad  \tx r_i=r_i\meas\Gamma
\,.\end{arr}\]
Here the mass $m$ of the particle is an objective scalar
and the velocity $\tx v=\tx\velo_\Gamma$
and on the right-hand side $\tx r$ satisfies the transformation rule
\[
  r_i\circ Y=\sum_{\bar i,\bar j\geq0}mY_{i\p\bar i\bar j}\,\tx v^*_{\bar i}\tx v^*_{\bar j}
  +\sum_{\bar i\geq0}Y_{i\p\bar i}\,r^*_{\bar i}
\,.\]
With this matrix distribution $T$ and this vector distribution $\tx r$
the system reads
for covariant vector valued test functions $\zeta$
\[\begin{arr}l
  0=\dprod{\zeta}{-\tx\div\,T+\tx r}{\dis(\UU)}
  =\dprod{\tx\D\zeta}{T}{\dis(\UU)}+\dprod{\zeta}{\tx r}{\dis(\UU)}
\\
  =\int_{\RR^4}\big( 
    \tx\D\zeta\ddd(m\tx v\up{\tx v}T) + \zeta\dd r
  \big)\d\meas{\Gamma}
\\
  =\int_{\RR^4}\Big(
    \sum_{i,j\geq0}\de_j\zeta_im\tx\velo_{\Gamma i}\tx\velo_{\Gamma j}
    + \sum_{i\geq0}\zeta_ir_i
  \Big)\d\meas{\Gamma}
\\
  =\int_{\Gamma}\sum_{i\geq0}\zeta_i\Big(
    -\div^\Gamma
    \Big(m\tx\velo_{\Gamma i}\frac{\tx\velo_\Gamma}{|\tx\velo_\Gamma|}\Big)
    +\frac{r_i}{|\tx\velo_\Gamma|}
  \Big)\d\haus1
\,.\end{arr}\]
For the second integral the regularity of $\Gamma$ and $m$ are required,
and for the last identity $\tx\velo_\Gamma\in\TT(\Gamma)$ is used.
\end{stmt}
\begin{prf}{}
The test function $\zeta$
is a covariant vector, that is $\zeta^*=\up{\tx\D{Y}}T\zeta\circ Y$,
hence
\[
  \zeta^*_{\bar i\p\bar j}=
  \sum_{i,j}Y_{i\p\bar i}Y_{j\p\bar j}\,\zeta_{i\p j}\circ Y
  +\sum_iY_{i\p\bar i\bar j}\,\zeta_i\circ Y
\,.\]
Then $g:=\tx\D\zeta\ddd T+\zeta\dd\tx r$ satisfies
\[\begin{arr}l
  g^*=\tx\D\zeta^*\ddd T^*+\zeta^*\dd\rate^*
  =\sum_{\bar i,\bar j}\zeta^*_{\bar i\p\bar j}T^*_{\bar i\bar j}
  +\sum_{\bar i}\zeta^*_{\bar i}r^*_{\bar i}
\\
  =\sum_{i,j}\zeta_{i\p j}\circ Y 
  \sum_{\bar i,\bar j}Y_{i\p\bar i}Y_{j\p\bar j}T^*_{\bar i\bar j}
  +\sum_i\zeta_i\circ Y\sum_{\bar i\bar j}
  Y_{i\p\bar i,\bar j}T^*_{\bar i\bar j}
\\\hfill
  +\sum_i\zeta_i\circ Y\sum_{\bar i}Y_{i\p\bar i}r^*_{\bar i}
\\
  =\sum_{i,j}\zeta_{i\p j}\circ Y\,T_{ij}\circ Y
  +\sum_i\zeta_i\circ Y\,r_{i}\circ Y
\\
  =(\tx\D\zeta\ddd T+\zeta\dd\rate)\circ Y
  =g\circ Y
\,,\end{arr}\]
that is, $g$ is an objective scalar.
Hence by \EQU{*.int}
\[
  \dprod{\zeta}{-\div\, T+\tx r}{\dis(\UU)}
  =\int_\Gamma g\d\meas{\Gamma}
  =\int_{\Gamma^*}g^*\d\meas{\Gamma^*}
  =\dprod{\zeta^*}{-\div\, T^*+\tx r^*}{\dis(\UU^*)}
\,.\]
Altogether this shows that we deal with the 4-momentum equation
(see also section~\REF{approximation..}).
\end{prf}

We remark that the mass equation is a special case of the 4-momentum equation
as shown in \REF{impuls.mass.}.

\begin{stmt}{compare}{Example}
Let the coordinates be $y=(t,x)\in\RR^4$ and let the moving point be given by
$\Gamma:=\set{(t,\xi(t))}{t\in\RR}$ and define $\tx\xi(t):=(t,\xi(t))$.
Then
$\tx{\dot\xi}=(1,\dot\xi)\in\TT(\Gamma)$,
and we assume that the derivative points in the same direction as $\velo_\Gamma$,
that is, for $y=\tx\xi(t)$  \[
  \tx\velo_\Gamma(y)  =\lambda(y)\tde{t}\tx\xi(t)
  ,\quad \lambda(y)>0
\,.\]
Then the 4-momentum equation for the evolving point $\Gamma$ is with $\tx v=\tx\velo_\Gamma$ for test functions $\zeta$
\[\begin{arr}l
  0=\dprod{\zeta}{-\div(m\tx v\up{\tx v}T)+\tx r}{\dis(\UU)}
  =m\tx v\up{\tx v}T\int_{\RR^4}\big( 
    \tx\D\zeta\ddd(m\tx v\up{\tx v}T) + \zeta\dd r
  \big)\d\meas{\Gamma}
\\
  =\int_{\RR}\big(
    \sum_{ij}\de_j\zeta_im\tx v_i\tx v_j+\sum_{i}\zeta_ir_i
  \big)(t,\xi(t))\frac{\d{t}}{\lambda(t,\xi(t))}
\\
  =\int_{\RR}\sum_{i}\Big(
    \tde{t}\big(\zeta_i(t,\xi(t))\big)\,m\tx v_i(t,\xi(t))
    +\zeta_i\frac{r_i}{\lambda}(t,\xi(t))
  \Big)\d{t}
\\
  =\int_{\RR}\sum_{i}\zeta_i\Big(
    -{\lambda}(t,\xi(t))\tde{t}(m\tx v_i(t,\xi(t)))+r_i(t,\xi(t))
  \Big)\frac{\d{t}}{\lambda(t,\xi(t))}
\\
  =\int_{\RR}\zeta\dd\Big(      -\lambda\tde{t}\Big(m\lambda\tde{t}\tx\xi(t)\Big)+r
  \Big)\frac{\d{t}}{\lambda}
\,.\end{arr}\]
This is with respect to an arbitrary $\ewelt\maps\RR^4\to\RR^4$.
\end{stmt}
\begin{prf}{}
Since $\Gamma\subset\RR^4$ is one-dimensional we have
$\tx v=\tx\velo_\Gamma=\lambda(1,\dot\xi)$, and $\lambda>0$ by assumption,
therefore
\[
  \meas\Gamma=\frac{\haus1\cutted\Gamma}{|\tx\velo_\Gamma|}
  =\frac{\haus1\cutted\Gamma}{\lambda\sqrt{1+|\dot\xi|^2}}
,\]
hence for every function $g$
\[
  \int g\d\meas\Gamma=\int_{\RR}\frac{g(t,\xi(t))}{\lambda(t,\xi(t))}\d{t}
\,.\]
Since
\[
  \tde{t}\big(\zeta_i(t,\xi(t))\big)=\tx\grad\zeta_i\dd(1,\dot\xi)
  =\frac1\lambda\tx\grad\zeta_i\dd\tx v
  =\frac1\lambda\sum_j\de_j\zeta_i\cdot\tx v_j
\,,\]
we obtain the result.
\end{prf}
This shows that the 4-momentum equation for a moving point
$t\mapsto\xi(t)$ is equivalent to the ODE
\begin{Equ}{ODE}
  \lambda\tde{t}\Big(m\lambda\tde{t}\tx\xi(t)\Big)=\tx r
\,,\end{Equ}
where $\lambda>0$ is given by $\lambda\tx{\dot\xi}=\tx\velo_\Gamma$.

\begin{stmt}{case1}{Case 1} If $\ewelt=\ee_0$, that is, $t$ is the
``normal'' time variable, because $\lambda=1$ and the ODE reads
\[
  \tde{t}\Big(m\tde{t}\tx\xi(t)\Big)=\tx r
\,.\]
Writing mass and momentum equations separately we get
\[
  \tde{t}m=\rate,\quad \tde{t}\Big(m\tde{t}\xi\Big)=\force,\quad
  \tx r=\begin{mat}\rate\\\force\end{mat}
.\]
\end{stmt}
\begin{prf}{}
It is $1=\ewelt\dd\velo_\Gamma=\lambda\ee_0\dd(1,\dot\xi)=\lambda$ and therefore
\[
  \begin{mat}\rate\\\force\end{mat}
  :=\tx r=\tde{t}\Big(m\tde{t}\tx\xi(t)\Big)
  =\tde{t}\Big(m\begin{mat}1\\\dot\xi\end{mat}\Big)
  =\begin{mat}\dot m\\\Dot{m\dot \xi}\end{mat}
.\]
\end{prf}

\begin{stmt}{case2}{Case 2} If $\ewelt$ is arbitrary and $\Gamma$ is an evolving curve, that is $\lambda>0$,
then $t$ is a time variable which is specific to the observer. 
\begin{enum}
\num{s}
If the factor $\lambda$ is known, the observer can
choose to a variable $s$ defined by $t=\trf{t}(s)$ and
\[
  \trf{t}'(s)=\lambda(\tx\xi(\trf{t}(s)))
\,.\]
Then $\Gamma=\set{\tx\xi(\trf{t}(s))}{s\in\RR}$ and
we obtain the differential equation
\[
  \tde{s}\Big(m\tde{s}\tx\xi(\trf{t}(s))\Big)=\tx r
\,.\]
\num{xi}
If as in \REF{time.unique.}
\[
  e_0=\begin{mat}\gamma\\\gamma V\end{mat}
  \,,\quad \ewelt=e_0'=\begin{mat}\gamma\\-\frac{\gamma}{\cc^2}V\end{mat}
  \,,\quad |V|<\cc
\,,\]
then $\dot\xi\dd V<\cc^2$, which is satisfied if $|\dot\xi|<\cc$.
\num{ode}
The ODE cannot be split as in \REF{case1.}.
But we have, if $\ewelt$ is a constant,
\[
  \tx r=\rate e_0 + \force
\]
according to \REF{impuls.reduct.}, a formula which is well known in
classical physics.
\end{enum}
\end{stmt}
\begin{prf}{\REF{s}}
For every function $t\mapsto h(t)$
\[
  \tde{s}h(\trf{t}(s))=\trf{t}'(s)\tde{t}h(t)=\lambda(\tx\xi)\tde{t}h(t)
\,.\]
\end{prf}
\begin{prf}{\REF{xi}} It is
\[
  1=\ewelt\dd\tx\velo_\Gamma=\lambda\,\ewelt\dd\tx{\dot\xi}
  =\begin{mat}\gamma\\-\frac{\gamma}{\cc^2}V\end{mat}\dd
  \begin{mat}1\\\dot\xi\end{mat}
  =\lambda\gamma\big(1-\frac1{\cc^2}V\dd\dot\xi\big)
\]
hence
\[
  0<\frac1{\lambda\gamma}=\big(1-\frac1{\cc^2}V\dd\dot\xi\big)
\,.\]
\end{prf}
\begin{prf}{\REF{ode}} The pure force $\force:=\tx r-\tx r\dd e_0'e_0$
has been defined in \REF{impuls.reduct.}. And $\tx r\dd e_0'=\ewelt$
by \REF{impuls.mass.} if $\ewelt$ is constant.
\end{prf}

The following example is based on the coordinates in \REF{case2.}.
You will find this in
\CIT{Einstein1905}{I {\S}5 Additionstheorem der Geschwindigkeiten}.

\begin{stmt}{addition}{Addition of velocities}
Let $y^*=(t^*,x^*)\in\RR^4$ and $y=(t,x)\in\RR^4$ be connected by a
Lorentz transformation $y=Y(y^*):={\bf L}_\cc(V,Q)y^*$
and consider two moving points $\Gamma:=\set{(t,\xi(t))}{t\in\RR}$
and $\Gamma^*:=\set{(t^*,\xi^*(t^*))}{t^*\in\RR}$
with $\Gamma=Y(\Gamma^*)$. Define
$u(t,\xi(t)):=\dot\xi(t)$ and $u^*(t^*,\xi^*(t^*)):=\dot\xi^*(t^*)$.
\begin{enum}
\num{o} Then 
\[
  u=\frac{\displaystyle V+\frac1{\gamma}{\bf B}_\cc(V)Qu^*}
  {\displaystyle 1+\frac1{\cc^2}V\dd Qu^*}
\]
if the denominator is positive.
\num{e} If $Q=\Id$ and $V\in\fcn{span}\{u^*\}$ then \[
  u=\frac{\displaystyle V+u^*}
  {\displaystyle 1+\frac1{\cc^2}V\dd u^*}
\]
if the denominator is positive.
\end{enum}

\end{stmt}
\begin{prf}{\REF{o}} The identity $(t,\xi(t))=Y(t^*,\xi^*(t^*))$ is
\[\begin{arr}l
  \begin{mat}t\\\xi(t)\end{mat}
  =\begin{mat}\gamma t^*+\frac\gamma{\cc^2}\up{V}TQ\xi^*(t^*)
  \\\gamma t^*V+{\bf B}_\cc(V)Q\xi^*(t^*)\end{mat}
\,.\end{arr}\]
If $s$ is a common parameter,
that is $(t(s),\xi(t(s)))=Y(t^*(s),\xi^*(t^*(s)))$,
then the derivative with respect to $s$ gives the result.
\end{prf}
\begin{prf}{\REF{e}} It is
$u^*=(u^*\dd\hat V)\,\hat V$ and ${\bf B}_\cc(V)\hat V=\gamma\hat V$.
\end{prf}

If one takes instead the velocities $\velo_\Gamma$ and  $\velo_{\Gamma^*}$
of the moving points, one has at corresponding points
$(t,\xi(t))$ and $(t^*,\xi^*(t^*))$ the relation
$\velo_\Gamma={\bf L}_\cc(V,Q)\velo_{\Gamma^*}$ and 
\[
  \velo_\Gamma=\frac{\tx u}{\ewelt\dd\tx u}
  \,,\quad \tx u:=\tx{\dot\xi}
  \,,\quad \tx\xi(t):=(t,\xi(t))
  \,,\quad \tx u=(1,u)
\,,\]
similarly for $\velo_{\Gamma^*}$. Compare the identity in \REF{scalar.u.u}.

{}
{}
\sect{approximation}{Approximation \DE{}\EN{of particles}}

A particle is a mass concentrated at an evolving point $\Gamma\subset\RR^4$.
The aim is to approximate this particle by a mass density in a neighbourhood of this curve
which converges to the particle mass on $\Gamma$.
This is necessary in order to prove that the distributional
equations of section \REF{point..}
are really the 4-momentum equations.

\medskip
In the general situation this means
that we have functions $g_\eps\maps\RR^4\to\RR$ with support in a
neighbourhood of $\Gamma$
such that for a limit function $g_\Gamma\maps\Gamma\to\RR$
\[
  g_\eps\leb4\to g_\Gamma\,\meas\Gamma
  \quad\tm{as $\eps\to0$, where}\quad
  \meas\Gamma:=\frac1{|\tx\velo_{\Gamma}|}\haus1\cutted\Gamma
\]
is the measure from \EQU{point.*.meas}. This means that
for test functions $\eta\in C_0^\infty(\RR^4)$
\begin{Equ}{limit}
  \int_{\RR^4} \eta g_\eps\d\leb4
\lto
  \int_\Gamma \eta g_\Gamma\d\meas\Gamma
  \quad\tm{as $\eps\to0$.}
\end{Equ}
In this situation we prove the following Theorem \REF{limit.}
for objective scalars $g_\eps$ in spacetime
which as $\eps\to0$ behave like a Dirac sequence around $\Gamma$.
Before we do so it is useful to prove the following quite general statements.

\begin{stmt}{e}{Lemma}
The vector $\ewelt=e'_0$ satisfies
\begin{enum}
\num{id}
$\displaystyle
  |e_0\wedge\cdots\wedge e_3|=\frac{|e_1\wedge\cdots\wedge e_3|}{|e'_0|}
$.
\num{scalar}
$|e_0\wedge\cdots\wedge e_3|$ is an objective scalar.

\begin{rem}{Remark} It is
$|e_0\wedge\cdots\wedge e_3|=|\ee_0\wedge\cdots\wedge\ee_3|=1$, if the observer
is connected with the standard Lorentz observer.
\end{rem}
\end{enum}
\end{stmt}
\begin{prf}{\REF{id}}
$\{e'_0\til e_3\}$ is a basis of $\RR^4$ and hence
$e_0=\mu e'_0 + \sum_{j=1}^n\nu_je_j$ and
$1=e'_0\dd e_0=\mu|e'_0|^2$. Therefore the vector
$\bar e:=e_0-\sum_{j=1}^n\nu_je_j=\mu e'_0$ satisfies
$\bar e\dd e_i=0$ for $i\geq1$. We conclude
\[
  |e_0\wedge\cdots\wedge e_3|=|\bar e\wedge e_1\wedge\cdots\wedge e_3|
  =|\bar e|\cdot|e_1\wedge\cdots\wedge e_3|
\]
where $|\bar e|=|\mu e'_0|=\up{|e'_0|}{-1}$
(see the proof of \REF{time.dual.}) 
\end{prf}
\begin{prf}{\REF{scalar}}
By \REF{change.basis-trans.}
we have the transformation rule $e_k\circ Y=\tx\D{Y}e_k^*$. Hence
\[\begin{arr}l
  |e_0\wedge\cdots\wedge e_3|\circ Y=|\de_{e_0^*}Y\wedge\cdots\wedge\de_{e_3^*}Y|
\\
  =|e_0^*\wedge\cdots\wedge e_3^*|\cdot|\de_0Y\wedge\cdots\wedge\de_3Y|
\end{arr}\]
and $|\de_0Y\wedge\cdots\wedge\de_3Y|=|\fcn{det}\tx\D{Y}|=1$.
\end{prf}

\begin{stmt}{W}{Lemma} For any function $h\maps\set{(y,z)}{z\in\welt(y)}\to\RR$
\[\begin{arr}l
  \int_{\welt(y)}h(y,z)\frac{\d\haus3(z)}
  {|e_1(y)\wedge\cdots\wedge e_3(y)|}
\\
  =\int_{\welt^*(y^*)}h(Y(y^*),\D{Y}z^*)\frac{\d\haus3(z^*)}
  {|e^*_1(y^*)\wedge\cdots\wedge e^*_3(y^*)|}
\,.\end{arr}\]
\end{stmt}
This lemma is applied to the case that $h^*(y^*,z^*):=h(Y(y^*),\D{Y}z^*)$.
\begin{prf}{}
It follows with the transformation $\D{Y}(y^*)\maps\welt^*(y^*)\to\welt(y)$
\[\begin{arr}c
  \int_{\welt(y)}h(y,z)\d\haus3(z)
  =|\fcn{det}\restrict{\D{Y}(y^*)}{\welt^*(y^*)}{}|
  \int_{\welt^*(y^*)}h(Y(y^*),\D{Y}z^*)\d\haus3(z^*)
\,.\end{arr}\]
If $\{e^\perp_1,e^\perp_2,e^\perp_3\}$ is an orthonormal basis of $\welt^*$ then
\[
  |\fcn{det}\restrict{\D{Y}}{\welt^*}{}|
  =|\de_{e^\perp_1}Y\wedge\cdots\wedge\de_{e^\perp_3}Y|
\]
Since $\{e^*_1,e^*_2,e^*_3\}$ is a basis of the same space $\welt^*$ we get
\[
  |\de_{e^*_1}Y\wedge\cdots\wedge\de_{e^*_3}Y|
  =|{e^*_1}\wedge\cdots\wedge{e^*_3}|\cdot
  |\de_{e^\perp_1}Y\wedge\cdots\wedge\de_{e^\perp_3}Y|
\,,\]
and since $\D{Y}e^*_i=e_i\circ Y$ for $i\geq1$
by \REF{change.basis-trans.} we see that
\[
  |\de_{e^*_1}Y\wedge\cdots\wedge\de_{e^*_3}Y|
  =|{e_1\circ Y}\wedge\cdots\wedge{e_3\circ Y}|
  =|{e_1}\wedge\cdots\wedge{e_3}|\circ Y
\,.\]
This gives
\[
  |\fcn{det}\restrict{\D{Y}}{\welt^*}{}|
  =\frac{|{e_1}\wedge\cdots\wedge{e_3}|\circ Y}
  {|{e^*_1}\wedge\cdots\wedge{e^*_3}|}
\]
and finishes the proof.
\end{prf}

We use this in order show that a mass point is the limit of an distributed
objective mass density.
\begin{stmt}{limit}{Theorem (\DE{}\EN{Convergence to} $\Gamma$)}
We assume that $\Gamma$ is an evolving point and let $g_\eps\maps\RR^4\to\RR$ be objective scalars whose support
is in the $\eps$-neigh\-bour\-hood of $\Gamma$.
We assume that for $y\in\Gamma$ and $z_\eps\to z\in\welt(y)$ \begin{equ}{conv}
  \eps^3|e_1(y)\wedge\cdots\wedge e_3(y)|\,g_\eps(y+\eps z_\eps)
  \to |\ewelt(y)|\,g(y,z) \quad\tm{as}\eps\to0
.\end{equ}
Then for $\eta\in C_0^\infty(\RR^4)$  \begin{Equ}{lim}\begin{arr}c
  \lim_{\eps\dto0}\int_{\RR^4}\eta g_\eps\d\leb4
  =\int_\Gamma \eta g_\Gamma\d\meas\Gamma  \,,\end{arr}\end{Equ}
where \begin{Equ}{gGamma}\begin{arr}c
  g_\Gamma(y)  :=\int_{\welt(y)}g(y,z)\frac{\d\haus3(z)}{|e_1(y)\wedge\cdots\wedge e_3(y)|}
\,.\end{arr}\end{Equ}
The function $g_\Gamma\maps\Gamma\to\RR$ is an objective scalar.
\end{stmt}
\begin{prf}{}
We let $\Gamma=\set{\tx\xi(s)}{s\in\RR}$ where the parameter $s$
can be chosen so that
$\tx\xi'(s)=\tx\velo_\Gamma(\xi(s))$. It follows that
for every function $h\maps\Gamma\to\RR$
\begin{Equ}{+}
  \int_\Gamma h(y)\d\meas\Gamma(y)
  =\int_\Gamma h(y)\frac{\d\haus1(y)}{|\tx\velo_\Gamma(y)|}
  =\int_\RR h(\tx\xi(s))\d\leb1(s)
\,.\end{Equ}
Now we use the transformation
\[
  (s,z)\mapsto y=\trf{y}(s,z):=\tx\xi(s)+\sum_{i\geq1}z_ie^\perp_i(\tx\xi(s))
  \in\RR^4
\,,\]
where $\{e^\perp_1(y),e^\perp_2(y),e^\perp_3(y)\}$
is an orthonormal basis of $\welt(y)$.
We compute its determinant,
since $\{\hat\ewelt,e^\perp_1,e^\perp_2,e^\perp_3\}$
is an orthonormal basis of $\RR^4$,
in an $\eps$-neighbourhood of $\Gamma$ as
\[\begin{arr}l
  |\fcn{det}\D\trf{y}(s,z)|
  =|(\tx\xi'+\sum_{i\geq1}z_i\D{e^\perp_i}(\tx\xi)\tx\xi')
  \wedge e^\perp_1\wedge\cdots\wedge e^\perp_3|
\\\hfill
  =|\tx\xi'\wedge e^\perp_1\wedge\cdots\wedge e^\perp_3|+\Order{\eps}
  =\tx\xi'\dd\hat\ewelt+\Order{\eps}
\,,\\
  \tx\xi'(s)\dd\hat\ewelt(\tx\xi(s))
  =\tx\velo_\Gamma(\xi(s))\dd\hat\ewelt(\tx\xi(s))
  =\frac1{|\ewelt(\tx\xi(s))|}
\,,\end{arr}\]
hence
\[
  |\fcn{det}\D\trf{y}(s,z)|
  =\frac1{|\ewelt(\tx\xi(s))|}+\Order{\eps}
.\]
It then follows that for $\eps\to0$
\[\begin{arr}l
  \int_{\RR^4}\eta g_\eps\d\leb4
  =\int_\RR\int_{\RR^3}(\eta g_\eps)(\trf{y}(s,z))|\fcn{det}\D\trf{y}(s,z)|
  \d\leb3(z)\d\leb1(s)
\\
  =\int_\RR\int_{\RR^3}(\eta g_\eps)(\trf{y}(s,z))
  \frac{\d\leb3(z)}{|\ewelt(\tx\xi(s))|}\d\leb1(s)+\Order{\eps}
\,.\end{arr}\]
Therefore we consider the following function
\[\begin{arr}l
  h_\eps(y):=\restrict{\int_{\RR^3}(\eta g_\eps)(\trf{y}(s,z))
  \frac{\d\leb3(z)}{|\ewelt(\tx\xi(s))|}}
  {\tx\xi(s)=y}{}
\\
  =\int_{\welt(y)}(\eta g_\eps)(y+z)\frac{\d\haus3(z)}{|\ewelt(y)|}
  =\int_{\welt(y)}\eps^3(\eta g_\eps)(y+\eps z)\frac{\d\haus3(z)}{|\ewelt(y)|}
\\
  =\eta(y)\int_{\welt(y)}\eps^3g_\eps(y+\eps z)\frac{\d\haus3(z)}{|\ewelt(y)|}
  +\Order{\eps}
\,,\end{arr}\]
and \EQU{+} gives
\[
  \int_{\RR^4}\eta g_\eps\d\leb4
  =\int_\Gamma h_\eps(y)\d\meas\Gamma(y)+\Order{\eps}
\,.\]
By assumption \EQU{conv} as $\eps\to0$ for $z\in\welt(y)$
\[
  \frac{\eps^3g_\eps(y+\eps z)}{|\ewelt(y)|}
  \to\frac{g(y,z)}{|e_1(y)\wedge\cdots\wedge e_3(y)|}
\,,\]
and this is why we obtain
\[
  h_\eps(y)\to\eta(y)
  \int_{\welt(y)}g(y,z)\frac{\d\haus3(z)}{|e_1(y)\wedge\cdots\wedge e_3(y)|}
  =\eta(y)g_\Gamma(y)=:h(y)
\,,\]
and this is the convergence as $\eps\to0$ of \EQU{lim}.

\medskip
It remains to prove that $g_\Gamma$ is an objective scalar. 
We know that $g_\eps$ is an objective scalar, hence $g_\eps(y)=g^*_\eps(y^*)$
for $y=Y(y^*)$.
By \REF{e.} the value
\[
  \lambda(y):=|e_0(y)\wedge\cdots\wedge e_3(y)|
  =\frac{|e_1(y)\wedge\cdots\wedge e_3(y)|}{|\ewelt(y)|}
\]
is an objective scalar, and therefore also
$\lambda(y)g_\eps(y)=\lambda^*(y^*)g^*_\eps(y^*)$
for $y=Y(y^*)$.
Therefore, if in addition $z=\D{Y}(y^*)z^*$, $z^*\in\welt^*(y^*)$,
from assumption \EQU{conv}
\[\begin{arr}l
  g^*(y^*,z^*)\ot\eps^3\lambda^*(y^*)g^*_\eps(y^*+\eps z^*)
  =\eps^3\lambda(y)g_\eps(Y(y^*+\eps z^*))
\\
  =\eps^3\lambda(y)g_\eps(y+\eps z_\eps) 
  \,,\quad z_\eps:=\frac1\eps(Y(y^*+\eps z^*)-Y(y^*))\,,
\\
  \to g(y,z)
  \,,\tm{since} z_\eps\to\D{Y}(y^*)z^*=:z\in\welt(y)\,,
\end{arr}\]
that is,
\begin{equ}{g}
  g^*(y^*,z^*)=g(y,z) \tm{for $y=Y(y^*)$, $z=\D{Y}(y^*)z^*$.}
\end{equ}
From this it follows $g_{\Gamma^*}(y^*)=g_{\Gamma}(y)$, since we prove in \REF{W.} that the integral which defines $g_{\Gamma}$ is
frame independent.
\end{prf}

We wanted to clarify the connection between the mass $m\maps\Gamma\to\RR$
of a moving point $\Gamma$ and
the mass density $\rho_\eps\maps\RR^4\to\RR$ in spacetime,
which is concentrated near $\Gamma$. We have shown in \REF{limit.} 
that a certain convergence of the usually called mass in the point $y\in\RR^4$
\[
  z\mapsto\eps^3|e_1(y)\wedge\cdots\wedge e_3(y)|\,\rho_\eps(y+\eps z)
\]
implies that $\rho_\eps\leb4\to m\meas\Gamma$ in distributional sense
where the mass of the particle in $y\in\Gamma$
\[
  m(y)=\frac1{|\ewelt(y)|}   \int_{\welt(y)}\lim_{\eps\to0}\big(\eps^3\rho_\eps(y+\eps z)\big)\d\haus3(z)
\]
is a mean value of the mass density across $\set{y+z}{z\in\welt(y)}$.
This mass $m$ is an objective scalar.
Accordingly, if the 4-velocity $\tx v_\eps$
converges strongly to $\tx\velo_\Gamma$ we have convergence of the mass
and the momentum equation.

{}}
\newpage
{}
\sect{fluid}{\DE{}\EN{Fluid equations}}
In this section we consider the 4-momentum equation \EQU{impuls.*.sys}
\[
  \tx\div\,T=\tau \quad\text{with tensor}\quad T=\rho\tx v\otimes\tx v+\tx\Pi
\,.\]
The characteristic behaviour of $T$ for a fluid is
that it depends on the gradients of the velocity $\tx\D\,\tx{v}$,
where here we take a 4-velocity $\tx{v}$ as defined in \REF{scalar.def-velo.}.
We prove the following theorem for the tensor $\tx\Pi$,
which is the generalization of the fact,
that in classical physics the dependence on the symmetric part of the velocity gradient is the only objective version of such a dependence.

\begin{stmt}{Dv}{Theorem} Let $\tx v$ be a 4-velocity
(as in \REF{scalar.def-velo.}).
\begin{enum}
\num{1}The tensor
$ \tilde S=
  \tx\D\tx v\,\GG+\up{(\tx\D\tx v\,\GG)}T-(\tx v\dd\tx\grad)\GG
$
is a contravariant tensor. \num{2}Also $ S:=
  \tx\D\tx v\,\GG^\raum+\up{(\tx\D\tx v\,\GG^\raum)}T-(\tx v\dd\tx\grad)\GG^\raum
$
is a contravariant tensor. \end{enum}
\end{stmt}

\begin{prf}{of contravariance}
Let $G$ be an arbitrary symmetric contravariant tensor,
that is a tensor with the property
\begin{equ}{o}
  G\circ Y=\tx\D{Y}G^*\up{\tx\D{Y}}T
\end{equ}
for observer transformations $y=Y(y^*)$.
The 4-velocity $\tx v$ satisfies 
\[
  \tx v_i\circ Y=\sum_{\bar i}Y_{i\p\bar i}\,\tx v^*_{\bar i}
\,.\]
Therefore one obtains for the derivative
\[
  \de_{\bar j}(\tx v_i\circ Y)
  =\sum_{\bar i}Y_{i\p\bar i\bar j}\,\tx v^*_{\bar i}
  +\sum_{\bar i}Y_{i\p\bar i}\de_{\bar j}\tx v^*_{\bar i}
\]
and from the chain rule
$
  \de_{\bar j}(\tx v_i\circ Y)
  =\sum_{j}(\de_j\tx v_i)\circ Y\,Y_{j\p\bar j}
$,
that is
\[
  \sum_{j}(\de_j\tx v_i)\circ Y\,Y_{j\p\bar j}
  =\sum_{\bar i}Y_{i\p\bar i\bar j}\,\tx v^*_{\bar i}
  +\sum_{\bar i}Y_{i\p\bar i}\de_{\bar j}\tx v^*_{\bar i}
\]
or in matrix notation
\[
  (\tx\D{v}\circ Y)\tx\D{Y}
  =\sum_{\bar i}\tx v^*_{\bar i}\tx\D{Y}_{\p\bar i}
  +\tx\D{Y}\tx\D{\tx v^*}
\,.\]
Multiplying this identity from the right side by $G^*\up{\tx\D{Y}}T$
one obtains using the property \EQU{o}
\[
  (\tx\D{v}\,G)\circ Y
  =\sum_{\bar i}\tx v^*_{\bar i}\tx\D{Y}_{\p\bar i}G^*\up{\tx\D{Y}}T
  +\tx\D{Y}(\tx\D{\tx v^*}G^*)\up{\tx\D{Y}}T
\,.\]
This is the transformation rule for $\tx\D{v}\,G$
(and it is, for $G=\GG_\cc$ and $c\to\infty$,
identical with the classical formula).
From this we obtain the transposed version ($G^*$ is symmetric)
\[
  \up{(\tx\D{v}\,G)}T\circ Y
  =\sum_{\bar i}\tx v^*_{\bar i}\tx\D{Y}G^*\up{\tx\D{Y}}T_{\p\bar i}
  +\tx\D{Y}\up{(\tx\D{\tx v^*}G^*)}T\up{\tx\D{Y}}T
\,.\]
The sum of both equations has as inhomogeneous term
\[
  M:=
  \sum_{\bar i}\tx v^*_{\bar i}\tx\D{Y}_{\p\bar i}G^*\up{\tx\D{Y}}T
  +\sum_{\bar i}\tx v^*_{\bar i}\tx\D{Y}G^*\up{\tx\D{Y}_{\p\bar i}}T
\,,\]
and reads
\[
  (\tx\D{v}\,G)\circ Y + \up{(\tx\D{v}\,G)}T\circ Y
  = M
  +\tx\D{Y}\big(\tx\D{\tx v^*}G^* + \up{(\tx\D{\tx v^*}G^*)}T\big)\up{\tx\D{Y}}T
\,.\]
The $M$-term also occurs in the transformation rule of $\sum_i\tx v_iG_{\p i}$,
since
\[\begin{arr}l
  \Big(\sum_i\tx v_iG_{\p i}\Big)\circ Y
  =\sum_{i\bar i}Y_{i\p\bar i}\tx v^*_{\bar i}(G_{\p i}\circ Y)
\\
  =\sum_{\bar i}\tx v^*_{\bar i}\Big(\sum_{i}Y_{i\p\bar i}G_{\p i}\circ Y\Big)
  =\sum_{\bar i}\tx v^*_{\bar i}(G\circ Y)_{\p\bar i}
  =\sum_{\bar i}\tx v^*_{\bar i}(\tx\D{Y}G^*\up{\tx\D{Y}}T)_{\p\bar i}
\\
  =\tx\D{Y}\Big(\sum_{\bar i}\tx v^*_{\bar i}G^*_{\p\bar i}\Big)\up{\tx\D{Y}}T
  +\Ubrack{\sum_{\bar i}\tx v^*_{\bar i}\Big(
  \tx\D{Y}_{\p\bar i}G^*\up{\tx\D{Y}}T
  +\tx\D{Y}G^*\up{\tx\D{Y}}T_{\p\bar i}
  \Big)}{=M}
\,.\end{arr}\]
Hence subtracting both equations gives $S:=\tx\D{v}\,G+\up{(\tx\D{v}\,G)}T-(\tx v\dd\grad)G$ and
this matrix satisfies $S\circ Y=\tx\D{Y}S^*\up{\tx\D{Y}}T$.
\end{prf}
\begin{prf}{\REF{1}}
Because $G:=\GG$ satisfies \EQU{o}.
\end{prf}
\begin{prf}{\REF{2}}
Also $G:=\GG^\raum=\GG-\GG^\zeit$ satisfies \EQU{o},
since $\GG^\zeit=-\frac1{\cc^2}e_0\up{e_0}T$ and $e_0$ is a contravariant vector.
\end{prf}

\begin{stmt}{J}{Theorem} Let $\tx v$ be a 4-velocity
as in \REF{scalar.def-velo.}.
\begin{enum}
\num{3}Define $\JJ:=\up{\ewelt}TS$ for the tensor $S$
in \REF{Dv.2}. The vector satisfies
\[
  \JJ_l=
  \sum_{ik}\ewelt_{i\p k}(\tx v_k\GG^\raum_{il}-\tx v_i\GG^\raum_{kl})
\]
and is a contravariant vector. \num{4}It is $\ewelt\dd\JJ=0$.
\end{enum}
\end{stmt}

\begin{prf}{\REF{3}}It is 
\[\begin{arr}l
  (\up\ewelt{T}S)_j=\sum_i\ewelt_iS_{ij}
\\
  =\sum_{ik}\ewelt_i\tx{v}_{i\p k}\GG^\raum_{kj}
  +\sum_k\Ubrack{\sum_i\ewelt_i\GG^\raum_{ik}}{=0}\tx{v}_{j\p k}
  -\sum_{ik}\tx{v}_{k}\ewelt_i\GG^\raum_{ij\p k}
\\
  =\sum_k\de_k\Ubrack{\big(\sum_i\ewelt_i\tx{v}_i\big)}{=1}\cdot\GG^\raum_{kj}
  -\sum_{ki}\ewelt_{i\p k}\tx{v}_i\GG^\raum_{kj}
\\
  -\sum_k\tx{v}_k\de_k\Ubrack{\big(\sum_i\ewelt_i\GG^\raum_{ij}\big)}{=0}
  +\sum_{ki}\tx{v}_k\ewelt_{i\p k}\GG^\raum_{ij}
\\
  =\sum_{ki}\ewelt_{i\p k}(\tx{v}_k\GG^\raum_{ij}-\tx{v}_i\GG^\raum_{kj})
\,.\end{arr}\]
Thus
\begin{equ}{3}
  (\up\ewelt{T}S)_j
  =\sum_{ki}\ewelt_{i\p k}(\tx{v}_k\GG^\raum_{ij}-\tx{v}_i\GG^\raum_{kj})
\,.\end{equ}
Now
\[
  \ewelt_{i\p k}=\frac{\ewelt_{i\p k}+\ewelt_{k\p i}}2+E_{ik}
\,,\quad
  E_{ik}=\frac{\ewelt_{i\p k}-\ewelt_{k\p i}}2
\,,\]
and since the bracket in \EQU{3} is antisymmetric in $(i,k)$, it follows that
\[
  \JJ_j:=(\up\ewelt{T}S)_j
  =\sum_{ki}E_{ik}(\tx{v}_k\GG^\raum_{ij}-\tx{v}_i\GG^\raum_{kj})
\,,\]
a term which we handled already in \REF{scalar.relativistic.}.
This is because $\ewelt^*_{\bar i}=\sum_iY_{i\p\bar i}\ewelt_i\circ Y$, and hence
\[
  \ewelt^*_{\bar i\p\bar k}=\sum_iY_{i\p\bar i\bar k}\ewelt_i\circ Y
  +\sum_{ik}Y_{i\p\bar i}Y_{k\p\bar k}\ewelt_{i\p k}\circ Y
\,.\]
Since $Y_{i\p\bar i\bar k}$ is symmetric in $(\bar i,\bar k)$, we obtain
\[
  E^*_{\bar i\bar k}=\sum_{ik}Y_{i\p\bar i}Y_{k\p\bar k}E_{ik}\circ Y
\,,\]
a property which was assumed in \REF{scalar.relativistic.}.
\end{prf}
\begin{prf}{\REF{4}}
Since $\sum_j\ewelt_j\GG^\raum_{kj}=0$ for every $k$.
\end{prf}

For fluids one has the following momentum system on the basis of \REF{Dv.2}
\[\begin{arr}c
  \tx\div(\rho\tx v\up{\tx v}T+\tx\Pi)=\tau
\,,\quad
  \tx\Pi=p\GG^\raum-S
\,,\\
  S=\mu
  (\tx\D\tx v\,\GG^\raum+\up{(\tx\D\tx v\,\GG^\raum)}T-(\tx v\dd\tx\grad)\GG^\raum)
  +\lambda\,\tx\div\tx v\,\GG^\raum
\,,\end{arr}\]
where $p$, $\mu$, and $\lambda$ are objective scalars.
It contains the mass equation
\[\begin{arr}c
  \tx\div(\rho\tx v+\JJ)=\rate
\,,\quad
  \JJ=\ewelt\dd\tx\Pi=\ewelt\dd S
\,,\end{arr}\]
where $\ewelt\dd\JJ=0$ and
$
  \rate=\ewelt\dd\tau+\D\ewelt\ddd\big(\tx v\up{(\rho\tx v+\JJ)}T+\tx\Pi\big)
$.

{}
\newpage
{}

\sect{higher}{\DE{}\EN{Higher moments}}
We now present the relativistic version of higher moments.
We consider moments of order $N$
with flux $T=\seq{T_\beta}{\beta\in\{0\til3\}^{N+1}}$.
Writing $T_\beta=T_{\alpha j}$ with
$\beta=(\alpha,j)$, where $\alpha\in\{0\til3\}^N$ and $j\in\{0\til3\}$,
the system of $N^{\rm th}$-moments reads in the version for test functions
$\zeta=\seq{\zeta_\alpha}{\alpha\in\{0\til3\}^N}$
\DE{}\EN{with} $\zeta_\alpha\in C_0^\infty(\UU;\RR)$
in a domain $\UU\subset\RR^4$
\begin{Equ}{sys}\begin{arr}l
  \int_{\RR^4}\sum_\alpha\Big(
  \sum_{j\geq0}\de_{y_j}\zeta_\alpha\cdot T_{\alpha j}
  +\zeta_\alpha\cdot g_\alpha
  \Big)\d\leb4=0
\,.\end{arr}\end{Equ}
The strong version of this system is $\tx\div\,T=g$ or
\begin{Equ}{sys-ind}
  \sum_{j\geq0}\de_{y_j}T_{\alpha j}=g_\alpha
  \quad\tm{\DE{}\EN{for}} \alpha\in\{0\til3\}^N
.\end{Equ}
\DE{}
\EN{The definition of the physical quantities of
this \DEF{system of $N^{\rm th}$-moments} is the following:
We demand from the test functions that they
satisfy the transformation rule}
\[\begin{arr}c
  \zeta^*_{\bar k_1\cdots\bar k_N}
  =\sum_{k_1\til k_N\geq0}
  Y_{k_1\p\bar k_1}\cdots Y_{k_N\p\bar k_N}
  \zeta_{k_1\cdots k_N}\circ Y
\end{arr}\]
for all $\bar k_1\til\bar k_N\in\{0\til3\}$,
\DE{}
\EN{where $Y$ is a relativistic observer transformation.}
\DE{}
\EN{Hence these test functions $\zeta$ are co\-va\-ri\-ant $N$-tensors.} 
Here $y$, $T$, $g$, and $\zeta$ are the quantities for one observer
and similarly $y^*$, $T^*$, $g^*$, and $\zeta^*$
are these quantities for another observer,
and $y=Y(y^*)$ is the observer transformation
where we as always assume that $\fcn{det}\tx\D{Y}=1$. Therefore it holds \[
  \int_{\RR^4}\sum_\alpha\Big(
  \sum_{j\geq0}\de_{y_j}\zeta_\alpha\cdot T_{\alpha j}
  +\zeta_\alpha\cdot g_\alpha
  \Big)\d\leb4
=
  \int_{\RR^4}\sum_\alpha\Big(
  \sum_{\bar j\geq0}\de_{y^*_{\bar j}}\zeta^*_\alpha\cdot T^*_{\alpha\bar j}
  +\zeta^*_\alpha\cdot g^*_\alpha
  \Big)\d\leb4
\]
and this is satisfied if the physical quantities $T$ and $g$ fulfill
the following transformation rule
(see the result in \REF{continuum.thm.} below)
\begin{Equ}{rule-T}
  T_{k_1\cdots k_M}\circ Y
  =\sum_{\bar k_1\til\bar k_M\geq0}
  Y_{k_1\p\bar k_1}\cdots Y_{k_M\p\bar k_M}
  T^*_{\bar k_1\cdots\bar k_M}
\end{Equ}
for $k_1\til k_M\in\{0\til3\}$ and $M=N+1$,
and
\begin{Equ}{rule-g}
  g_{k_1\cdots k_N}\circ Y
  =\sum_{\bar k_1\til\bar k_N,j\geq0}
  \big(Y_{k_1\p\bar k_1}\cdots Y_{k_N\p\bar k_N}\big)_{\p j}
  T^*_{\bar k_1\cdots\bar k_N j}
\\\hfill
  +\sum_{\bar k_1\til\bar k_N\geq0}
  Y_{k_1\p\bar k_1}\cdots Y_{k_N\p\bar k_N}
  g^*_{\bar k_1\cdots\bar k_N}
\end{Equ}
for $k_1\til k_N\in\{0\til3\}$.
If one wants to express a high moment in terms of another observer,
one needs all moments of the other observer up to this order.
As symmetry condition one might assume that $T_{\alpha j}$ and $g_\alpha$
are symmetric in the components of $\alpha$, but in general there is
no symmetry with respect to the last index $j$. Thus for $N=1$ the
relativistic Navier-Stokes equations are included. And it is important to say
that also for arbitrary $N$ we do not prescribe a constitutive relation
for $T_\beta$, we only assume that they are a solution of system \EQU{sys}.
The form of this system is the only connection to section \REF{chapman..}.

\medskip
So far no special relativistic argument has occurred.
But obviously the question arises, how the lower order momentum equations
are contained in this presentation, are the $(N-1)^{\rm th}$-moments
part of the $N^{\rm th}$-moments as in the non-relativistic case?
So we have to find a relativistic version
of the reduction in \REF{chapman.reduct.}.
On the other hand what is clear is that the usual representation of the tensor
\begin{equ}{vv}
  T_\beta=\rho{\tx v}_{\beta_1}\cdots{\tx v}_{\beta_M}+\tx\Pi_\beta
\end{equ}
with $M=N+1$ can be used also in the relativistic case.
Here we choose the 4-velocity $\tx v$ as defined in \REF{scalar.def-velo.}
and $\rho$ as the mass density which is an objective scalar.
Then the tensor $T$ satisfies \EQU{rule-T}, if $\tx\Pi$ does it,
because $\tx v$ is a contravariant vector, that is,
\[
  {\tx v}_i\circ Y=\sum_{\bar i=0}^3Y_{i\p\bar i}{\tx v}^*_{\bar i}
\,.\]
This implies
\[\begin{arr}l
  \Big(\rho\prod_{i=1}^M{\tx v}_{k_i}\Big)\circ Y
  =\rho\circ Y\prod_{i=1}^M{\tx v}_{k_i}\circ Y
  =\rho^*\prod_{i=1}^M\Big(
  \sum_{\bar k_i=0}^3Y_{k_i\p\bar k_i}{\tx v}^*_{\bar k_i}\Big)
\\
  =\rho^*\sum_{\bar k_1\til\bar k_M=0}^3
  Y_{k_1\p\bar k_1}\cdots Y_{k_M\p\bar k_M}
  {\tx v}^*_{\bar k_1}\cdots{\tx v}^*_{\bar k_M}
\\
  =\sum_{\bar k_1\til\bar k_M=0}^3
  Y_{k_1\p\bar k_1}\cdots Y_{k_M\p\bar k_M}
  \Big(\rho^*\prod_{i=1}^M\tx v^*_{\bar k_i}\Big)
\end{arr}\]
and gives \EQU{rule-T} for $T$.
 
\medskip
\TITEL{Reduction of the system}
The system of $N^{\rm th}$-moments \EQU{sys-ind}
consists of $4^N$ differential equations
(without taking symmetry into account)
and it should contain the $4^{N-1}$ differential equations
of the $(N-1)^{\rm th}$-moments equation.
We realize this the same way as we did for $N=1$ in section \REF{impuls..}.
There we considered the 4-moment system and we showed in \REF{impuls.reduct.}
that it contains the mass equation.
Here we present a generalization.
We use as special test function
\begin{equ}{eta}
  \zeta_{\alpha_1\cdots\alpha_N}:=\ewelt_{\alpha_1}\eta_{\alpha_2\cdots\alpha_N}
\end{equ}
where the vector $\ewelt$ is the time vector from section \REF{time..}.
The function $\eta$ is a covariant $(N-1)$-tensor.

\begin{stmt}{reduction}{Reduction lemma}
The system of $N^{\rm th}$-moments \EQU{*.sys-ind} contains as part
the system of $(N-1)^{\rm th}$-moments
\[
  \sum_{j\geq0}\de_{y_j}T_{\alpha j}=g_{\alpha}
  \quad\tm{\DE{}\EN{for}} \alpha\in\{0\til3\}^{N-1}
\]
with
\[\begin{arr}l
  T_{\alpha j}:=\sum_i\ewelt_iT_{i\alpha j}
\,,\quad
  g_\alpha:=\sum_i\big(\ewelt_{i\p j}\,T_{i\alpha j}+\ewelt_i\,g_{i\alpha}\big)
\,.\end{arr}\]
\begin{rem}{Remark}
The case $N=1$ is included by writing $\sum_{j\geq0}\de_{y_j}T_{j}=g$.
\end{rem}
\end{stmt}
Since $\sum_i\ewelt_i\tx v_i=1$, this is consistent with \EQU{*.vv},
that is, \EQU{*.vv} holds for all orders of moments.
\begin{prf}{} We choose the test function as in \EQU{*.eta}.
If $\eta$ is a covariant $(N-1)$-tensor
then since $\ewelt$ is a covariant vector
\[\begin{arr}l
  \zeta^*_{\bar\alpha_1\cdots\bar\alpha_N}
  =\ewelt^*_{\bar\alpha_1}\eta^*_{\bar\alpha_2\cdots\bar\alpha_N}
\\
  =\sum_{\alpha_1\geq0}Y_{\alpha_1\p\bar\alpha_1}\ewelt_{\alpha_1}\circ Y
  \cdot
  \sum_{\alpha_2\til\alpha_N\geq0}Y_{\alpha_2\p\bar\alpha_2}\cdots Y_{\alpha_N\p\bar\alpha_N}
  \eta_{\alpha_2\cdots\alpha_N}\circ Y
\\
  =\sum_{\alpha_1\til\alpha_N\geq0}Y_{\alpha_1\p\bar\alpha_1}\cdots Y_{\alpha_N\p\bar\alpha_N}
  \ewelt_{\alpha_1}\circ Y\eta_{\alpha_2\cdots\alpha_N}\circ Y
\\
  =\sum_{\alpha_1\til\alpha_N\geq0}Y_{\alpha_1\p\bar\alpha_1}\cdots Y_{\alpha_N\p\bar\alpha_N}
  \zeta_{\alpha_1\cdots\alpha_N}\circ Y
\,.\end{arr}\]
This means that $\zeta$ is an allowed test function and
it follows from \EQU{*.sys} writing $\alpha=(i,\gamma)$
\[\begin{arr}l
  0=
  \int_{\RR^4}\sum_\alpha\Big(
  \sum_j\de_{y_j}\zeta_\alpha\cdot T_{\alpha j}
  +\zeta_\alpha\cdot g_\alpha
  \Big)\d\leb4
\\
  =\int_{\RR^4}\sum_{i\gamma}\Big(
  \sum_j\de_{y_j}(\eta_\gamma\ewelt_i) T_{i\gamma j}+\eta_\gamma\ewelt_i\cdot g_{i\gamma}
  \Big)\d\leb4
\\
  =\int_{\RR^4}\sum_\gamma\Big(
  \sum_j\de_{y_j}\eta_\gamma\cdot\Ubrack{\sum_i\ewelt_iT_{i\gamma j}}{=:T_{\gamma j}}
  +\eta_\gamma\Ubrack{\sum_i\Big(
  \de_{y_j}\ewelt_iT_{i\gamma j}+\ewelt_ig_{i\gamma}\Big)
  }{=:g_\gamma}
  \Big)\d\leb4
\,.\end{arr}\]
This gives the result.
\end{prf}

\medskip
\TITEL{Coriolis coefficients}
The transformation formula \EQU{rule-g} 
gives rise to the following definition of the coefficients
$\coriolis_\alpha=\seq{\coriolis^\beta_\alpha}{\beta\in\{0\til3\}^{N+1}}$
(this is a generalization of \REF{impuls.coriolis.})
\begin{Equ}{coriolis}\begin{arr}l  g_\alpha=\tx\force_\alpha+\sum_{\beta\in\{0\til3\}^{N+1}}\coriolis^\beta_\alpha T_\beta
\quad  \tm{for}\alpha\in\{0\til3\}^N .\end{arr}\end{Equ}
With these Coriolis coefficients the system \EQU{sys-ind} has the form
\begin{Equ}{sys-coriolis}
  \sum_{j\geq0}\de_{y_j}T_{\alpha j}
  -\sum_{\beta\in\{0\til3\}^{N+1}}\coriolis^\beta_\alpha T_\beta=\tx\force_\alpha
  \quad\tm{for} \alpha\in\{0\til3\}^N
\end{Equ}
with transformation rule \EQU{rule-T} for the tensor $T$, that is,
$T$ is a contravariant $M$-tensor ($M=N+1$), and
\begin{Equ}{rule-f}
  \tx\force_{k_1\cdots k_N}\circ Y
  =\sum_{\bar k_1\til\bar k_N\geq0}
  Y_{k_1\p\bar k_1}\cdots Y_{k_N\p\bar k_N}
  \,\tx\force^*_{\bar k_1\cdots\bar k_N}
\end{Equ}
for the force (e.g.~the gravity or the Lorentz force),
that is, the entire force
$\tx\force:=\seq{\tx\force_\alpha}{\alpha\in\{0\til3\}^N}$
is a contravariant $N$-tensor.
(See the equation \CIT{Mueller1998}{Chap.2 (3.15)} for a comparison
with the classical case.)
The Coriolis coefficients satisfy the
following transformation rule.

\begin{stmt}{coriolis}{Rule for the Coriolis coefficients}
The rule \EQU{*.rule-g} for $g$ is equivalent to the fact that $\tx\force$
is a contravariant $N$-tensor and
\[\begin{arr}l
  \sum_{m_1\til m_{N+1}\geq0}
  Y_{m_1\p\bar m_1}\cdots Y_{m_{N+1}\p\bar m_{N+1}}
  \coriolis^{m_1\cdots m_{N+1}}_{k_1\cdots k_N}\circ Y
\\
  =\sum_{\bar k_1\til\bar k_N\geq0}
  Y_{k_1\p\bar k_1}\cdots Y_{k_N\p\bar k_N}
  \coriolis^{*\bar m_1\cdots\bar m_{N+1}}_{\bar k_1\cdots\bar k_N}
  +\big(Y_{k_1\p\bar m_1}\cdots Y_{k_N\p\bar m_N}\big)_{\p\bar m_{N+1}}
\end{arr}\]
for all $k_1\til k_N$ and $\bar m_1\til\bar m_{N+1}$.
\end{stmt}
\begin{prf}{}
We take the equation \EQU{*.rule-g}.
Using the definition \EQU{*.coriolis} and
the above transformation rule \EQU{*.rule-f} for $\tx\force$
this equation becomes for $k_1\til k_N\in\{0\til3\}$
\begin{Equ}{cpre}\begin{arr}l
  \Big(
  \sum_{m_1\til m_{N+1}\geq0}\coriolis^{m_1\cdots m_{N+1}}_{k_1\cdots k_N} T_{m_1\cdots m_{N+1}}
  \Big)\,\circ\,Y
\\
  =\sum_{\bar k_1\til\bar k_N,j\geq0}
  \big(Y_{k_1\p\bar k_1}\cdots Y_{k_N\p\bar k_N}\big)_{\p j}
  T^*_{\bar k_1\cdots\bar k_N j}
\\\hfill
  +\sum_{\bar k_1\til\bar k_N\geq0}
  Y_{k_1\p\bar k_1}\cdots Y_{k_N\p\bar k_N}
  \sum_{\bar m_1\til\bar m_{N+1}\geq0}
  \coriolis^{*\bar m_1\cdots\bar m_{N+1}}_{\bar k_1\cdots\bar k_N} T^*_{\bar m_1\cdots\bar m_{N+1}}
\end{arr}\end{Equ}
Using \EQU{*.rule-T}, that is
\[
  T_{m_1\cdots m_{N+1}}\circ Y =\sum_{\bar m_1\til\bar m_{N+1}\geq0}
  Y_{m_1\p\bar m_1}\cdots Y_{m_{N+1}\p\bar m_{N+1}}T^*_{\bar m_1\cdots\bar m_{N+1}}
\,,\]
the left-hand side of \EQU{cpre} becomes
\[
  \sum_{\begin{sarr}l\bar m_1\til\bar m_{N+1}\geq0,\\m_1\til m_{N+1}\geq0\end{sarr}}
  Y_{m_1\p\bar m_1}\cdots Y_{m_{N+1}\p\bar m_{N+1}}
  \coriolis^{m_1\cdots m_{N+1}}_{k_1\cdots k_N}\circ Y
  \,T^*_{\bar m_1\cdots\bar m_{N+1}}
\,.\]
Now compare the coefficients of $T^*$
with the one of the right-hand side of \EQU{cpre}
and obtain
\[\begin{arr}l
  \sum_{m_1\til m_{N+1}\geq0}
  Y_{m_1\p\bar m_1}\cdots Y_{m_{N+1}\p\bar m_{N+1}}
  \coriolis^{m_1\cdots m_{N+1}}_{k_1\cdots k_N}\circ Y
\\
  =\big(Y_{k_1\p\bar m_1}\cdots Y_{k_N\p\bar m_N}\big)_{\p\bar m_{N+1}}
  +\sum_{\bar k_1\til\bar k_N\geq0}
  Y_{k_1\p\bar k_1}\cdots Y_{k_N\p\bar k_N}
  \coriolis^{*\bar m_1\cdots\bar m_{N+1}}_{\bar k_1\cdots\bar k_N}
\,,\end{arr}\]
which is the assertion.
\end{prf}

{}
{}
\sect{continuum}{\DE{}\EN{Appendix: Divergence systems}}
\newcommand{\rates}{r}
\newcommand{\M}{M}

We consider a spacetime domain $\UU\subset\RR^{n+1}$, $n=3$,
and in $\UU$ integrable fluxes $q^k$ and functions $\rates^k$, $k=0\til\M$,
which solve the divergence system in $\UU$ \begin{Equ}{system}
  \sum_{i=0}^n \de_{y_i} q_i^k =\rates^k
  \tm{for} k=0\til\M
.\end{Equ}
Further, we suppose that an invertible matrix $Z=Z(y^*)$,
\begin{equ}{Z}
  Z(y^*)=\seq{Z_{kl}(y^*)}{k,l=0\til\M}
\,,\end{equ}
is given. We consider the following transformation rule
for observer transformations $y=Y(y^*)$
\begin{equ}{o}\begin{arr}l
  q_i^k\circ Y=\frac1J\sum_{jl}Y_{i\p j}Z_{kl}q^{*l}_j  
\,,\quad
  J:=\fcn{det}\D_{y^*}Y > 0   \,,\\
  \rates^k\circ Y=\frac1J\left(
  \sum_{jl}Z_{kl\p j}q^{*l}_j
  +\sum_{l}Z_{kl}\rates^{*l}
  \right)
\,,\\
  \text{for all~} i=0\til n \tm{and} k=0\til\M
\,,\\
  \text{where $j$ runs from $0$ to $n$, and $l$ from $0$ to $\M$}
.\end{arr}\end{equ}
In \CIT{Alt-Kontinuum}{Section I.5} it has been proved
that the system \EQU{system} is invariant under observer transformations
if \EQU{o} is satisfied:

\begin{stmt}{thm}{Theorem}
If the quantities $q^k$, $\rates^k$, $k=0\til\M$,
satisfy the transformation rule \EQU{*.o} for a matrix $Z$ as in \EQU{*.Z},
then with $\UU=Y(\UU^*)$
\begin{Equ}{o}
  \sum_{l=0}^\M\int_{\UU^*}\Big(
    \sum_{j=0}^n\de_{y^*_j}{\zeta^*_l}\,q^{*l}_j+{\zeta^*_l}\rates^{*l}
  \Big)\d\leb{n+1}
\\
  =\sum_{k=0}^\M\int_{\UU}\Big(
    \sum_{i=0}^n\de_{y_i}{\zeta_k}\,q_i^k+{\zeta_k}\rates^k
  \Big)\d\leb{n+1}
\end{Equ}
where the test functions satisfy
\begin{equ}{zeta}
  \zeta^*=\up{Z}{T}\zeta\circ Y  
\,.\end{equ}
\end{stmt}

In the case $Z=\tx\D{Y}$ this theorem is used in this paper
for the relativistic theory special for the 4-momentum system.
In this case the condition \EQU{o} on the fluxes are for $i,k=0\til n$
\[
  q_i^k\circ Y=\frac1J\sum_{j,l=0}^nY_{i\p j}Y_{k\p l}q^{*l}_j  
\,,\]
and the test function $\zeta$ is a covariant vector.
For the hierarchical theory the matrix $Z$ is
$
  Z=\seq{Z_{(i_1\til i_N)(\bar i_1\til\bar i_N)}}
  {i_1\til i_N,\bar i_1\til\bar i_N=0\til n}
$ with
\[
  Z_{(i_1\til i_N)(\bar i_1\til\bar i_N)}
  =
  Y_{i_1\p\bar i_1}\cdots Y_{i_N\p\bar i_N}
\,.\]
In this case the property \EQU{thm.zeta} says that the test function
is a covariant $N$-tensor.
In the special case $Z=\Id$ this theorem can be used
for the introduction to elasticity theory,
see e.g.~\CIT{Alt-Kontinuum}{Section I.6}.

{}
{}
\sect{invariance}{\DE{}\EN{Appendix: Theorem on Lorentz matrix}} 
The following is a well known theorem.

\begin{stmt}{thm}{Theorem}
The following sets of matrices are the same.
\begin{enum}
\num{2}
The set of all matrices $M$ satisfying
\[
  \GG_\cc = M \GG_\cc \up{M}{T}
\]
with the normalization, that $M_{00}\geq 0$ and $\fcn{det}M>0$.
\num{1}
The set of all matrices
\[
  M={\bf L}_\cc(V,Q) 
\]
with $V\in\RR^3$, $|V|<\cc$,
and $Q$ an orthonormal matrix with determinant $1$.
\end{enum}
\end{stmt}

The Lorentz matrices ${\bf L}_\cc(V,Q)$ are given by
\[
  {\bf L}_\cc(V,Q)=
  \begin{mat}
    \gamma
  \tab
    \frac{\gamma}{\cc^2}\up{V}{T}Q
  \\
    \gamma V
  \tab
    {\bf B}_\cc(V)Q    \end{mat}
,\]
\DE{}\EN{where}
$
  {\bf B}_\cc(V):=
  \Id+\frac{\gamma^2}{\cc^2(\gamma+1)}V\up{V}{T}
$
and
$
  \gamma=\big(1-\frac{|V|^2}{\cc^2}\big)^{-\frac12}
  \tm{for}|V|<\cc
$.
{}
{}
\newcommand{\authors}[1]{#1:}
\newcommand{\titles}[1]{\emph{#1}.}
\newcommand{\sources}[1]{#1.}
\newcommand{\published}[1]{ #1}
\newcommand{\electro}[1]{\\#1}

{}

\begin{center}
\hrulefill\\[2ex]
\end{center}
\timstamp

\end{document}